# PROPORTIONAL HAZARDS MODELS WITH CONTINUOUS MARKS

By Yanqing Sun[1], Peter B. Gilbert[2] and Ian W. McKeague[3]

*University of North Carolina at Charlotte, University of Washington and Fred Hutchinson Cancer Research Center, and Columbia University*

For time-to-event data with finitely many competing risks, the proportional hazards model has been a popular tool for relating the cause-specific outcomes to covariates [Prentice et al. *Biometrics* **34** (1978) 541–554]. This article studies an extension of this approach to allow a continuum of competing risks, in which the cause of failure is replaced by a *continuous* mark only observed at the failure time. We develop inference for the proportional hazards model in which the regression parameters depend nonparametrically on the mark and the baseline hazard depends nonparametrically on both time and mark. This work is motivated by the need to assess HIV vaccine efficacy, while taking into account the genetic divergence of infecting HIV viruses in trial participants from the HIV strain that is contained in the vaccine, and adjusting for covariate effects. Mark-specific vaccine efficacy is expressed in terms of one of the regression functions in the mark-specific proportional hazards model. The new approach is evaluated in simulations and applied to the first HIV vaccine efficacy trial.

**1. Introduction.** It has been 30 years since Prentice et al. [14] introduced a Cox regression framework for the analysis of failure time data in the presence of finitely many competing risks. Yet many important applications of competing risks methodology involve *continuous* causes-of-failure (marks). In HIV vaccine trials, for example, genetic divergence of infecting

Received December 2006; revised September 2007.
[1]Supported in part by NSF Grant DMS-0604576, NIH Grant 2 RO1 AI054165-04 and funds provided by the University of North Carolina at Charlotte.
[2]Supported in part by NIH Grant 2 RO1 AI054165-04.
[3]Supported in part by NSF Grant DMS-0505201.
*AMS 2000 subject classifications.* Primary 62N01; secondary 62N02, 62N03, 62G20.
*Key words and phrases.* Competing risks, distribution-free confidence bands and tests, failure time data, genetic data, HIV vaccine trial, pointwise and simultaneous confidence bands, semiparametric model, survival analysis.

HIV viruses from the HIV strain represented in the vaccine needs to be taken into account to properly assess vaccine efficacy, but the mark variable is essentially continuous because of the large number of mutations involved. Other examples of continuous mark variables include lifetime medical cost or a quality of life score associated with survival time [13]. The grouping of continuous mark data into discrete marks is unsatisfactory because that amounts to a coarsening of the data and the results will depend on the way the groups are defined. To address this problem, we develop inference for a proportional hazards model in which both the regression parameters and the baseline hazard function depend nonparametrically on a continuous mark.

The paper is motivated by the need for new methods to analyze data from HIV vaccine efficacy trials. Approximately 15,000 new HIV infections occur each day [20], making development of a protective HIV vaccine a top priority for biomedical science. In efficacy trials thousands of HIV-negative volunteers are randomized to receive vaccine or placebo, and are monitored for HIV infection. Five efficacy trials have recently been conducted. A primary objective of each trial is to assess vaccine efficacy (VE) to prevent infection, where typically VE is defined as one minus the hazard ratio (vaccine/placebo) of HIV infection. One of the greatest barriers to achieving an efficacious vaccine is the extreme genetic heterogeneity of HIV [6, 11]. Although it may be possible to develop a vaccine that protects against HIV strains genetically similar to the HIV virus or viruses represented in the vaccine, it may be quite difficult to develop one to protect against HIV strains dissimilar from the vaccine material. This phenomenon is well known for flu vaccines—moderate genetic mismatch between an exposing flu virus and the virus represented in the vaccine causes vaccine failure, which has necessitated development of a new vaccine each year that is closely matched to the contemporary circulating flu strains. The genetic divergence (or distance) between two aligned HIV sequences can be measured as the weighted percent mismatch of amino acids, and since this distance may be unique for all infected subjects, it is natural to consider it as a continuous mark variable. The formidable problem of HIV genetic diversity implies that an important objective of an efficacy trial is assessment of if and how VE depends on the genetic divergence.

This problem can be addressed in terms of the *conditional mark-specific hazard function*, defined as

$$\lambda(t,v|z) = \lim_{h_1,h_2 \to 0} P\{T \in [t, t+h_1),$$
(1)
$$V \in [v, v+h_2)|T \geq t, Z(t) = z\}/h_1 h_2,$$

where $T$ is the failure (infection diagnosis) time, $V$ is a continuous mark variable and $Z(t)$ is a (possibly time-dependent) $p$-dimensional covariate. Huang and Louis [7] developed the nonparametric maximum likelihood estimator of the joint distribution of $T$ and $V$ in terms of the unconditional mark-specific hazard function. Gilbert, McKeague and Sun [5] defined mark-specific vaccine efficacy as $\mathrm{VE}(t,v) = 1 - \lambda(t,v|1)/\lambda(t,v|0)$, with $z$ being the indicator of membership in the vaccine group; they developed several nonparametric and semiparametric tests concerning $\mathrm{VE}(t,v)$.

In this article, we develop the mark-specific proportional hazards (PH) model

$$(2) \qquad \lambda(t,v|z(t)) = \lambda_0(t,v)\exp\{\beta(v)^T z(t)\},$$

where the baseline hazard function $\lambda_0(\cdot,v)$ and the $p$-dimensional regression parameter $\beta(v)$ are unknown continuous functions of $v$. As far as we know, this model has never been studied in the literature, even though it is closely related to the discrete cause-of-failure models discussed by Prentice et al. [14]. The approach in the continuous case departs from the discrete case in that it is necessary to "borrow strength" from data in a neighborhood of $v$, with the data closest to $v$ contributing the most.

For the HIV vaccine trial application, we partition the covariate as $z(t) = (z_1, z_2(t))^T$, where $z_1$ is the treatment (vaccine) group indicator and $z_2(t)$ is a vector of possibly time-dependent covariates. Then the vaccine efficacy defined above takes the simpler form $\mathrm{VE}(v) = 1 - \exp(\beta_1(v))$, without any dependence on $t$. By assuming proportional hazards, model (2) can provide more powerful tests of mark-specific vaccine efficacy than the nonparametric procedures of Gilbert, McKeague and Sun [5], and the model allows adjustment for covariate effects. Furthermore, ignoring the mark variable and studying vaccine efficacy using the standard Cox model, as is widely practiced in vaccine trials for many infectious diseases, can give misleading results. In fact, even in the case of model (2) with $z$ as the treatment indicator, the ordinary (marginal) Cox model will be misspecified unless the baseline $\lambda_0(t,v)$ factors into separate functions of $t$ and $v$.

Indeed, consider the model $\lambda(t,v|z=0) = \gamma_0/2 + \gamma_1 tv$ and $\lambda(t,v|z=1) = \gamma_0 v + \gamma_1 tv^2$, for $t \geq 0$, $0 \leq v \leq 1$, $z \in \{0,1\}$. The corresponding marginal hazard functions are $\lambda(t|z=0) = \gamma_0/2 + \gamma_1 t/2$ and $\lambda(t|z=1) = \gamma_0/2 + \gamma_1 t/3$, for $t \geq 0$. It is clear that $\lambda(t|z)$ is not a proportional hazards model unless $\gamma_0$ or $\gamma_1$ is zero. If $\gamma_1 = 0$, the resulting marginal hazards become proportional for $z=0$ and $z=1$. However, in this example, the marginal vaccine efficacy $\mathrm{VE} = 1 - \lambda(t|z=1)/\lambda(t|z=0) = 0$ while the mark-specific vaccine efficacy is $\mathrm{VE}(v) = 1 - 2v$. The ordinary Cox model averages the mark-specific vaccine efficacy over its range, and important vaccine effects may be missed. This issue will be further illustrated in our simulation study. In general, use of the ordinary Cox model for studying hazard ratios can be misleading if an important mark variable is ignored. The mark-specific PH model offers a way to correct for that deficiency.

We also consider a cumulative vaccine efficacy estimand defined as $\mathrm{CV}(v) = \int_a^v \mathrm{VE}(u)\,du$ where $a > 0$. We develop distribution-free uniform confidence bands for $\mathrm{CV}(v)$, which are useful for inferential purposes. In addition we derive test statistics for evaluating mark-specific vaccine efficacy based on the estimator of $\mathrm{CV}(v)$.

The paper is organized as follows. Section 2 develops a local partial likelihood procedure for estimating $\beta(v)$, leading to the construction of pointwise confidence intervals and formal tests for various hypotheses of interest concerning vaccine efficacy. A simulation study evaluating the performance of the proposed tests and the pointwise and simultaneous confidence intervals for $\mathrm{VE}(v)$ and $\mathrm{CV}(v)$ is presented in Section 3. The proposed methods are applied to analyze the data from the first HIV vaccine efficacy trial in Section 4. We discuss some general aspects of mark-specific PH models in Section 5. Proofs of the main results are placed in the Appendix.

## 2. Mark-specific proportional hazards model.

2.1. *Local partial likelihood.* We begin by stating some assumptions and notations that are used throughout the paper. The mark variable $V$ is assumed to have a known and bounded support; rescaling $V$ if necessary, this support is taken without loss of generality to be $[0, 1]$. The observations $(X_i, \delta_i, \delta_i V_i, Z_i)$, $i = 1, \ldots, n$, are assumed to be i.i.d. replicates of $(X, \delta, \delta V, Z)$, where $X$ is the right-censored failure time corresponding to $T$, which satisfies the model (2), and $\delta$ is the indicator of non-censorship. The mark is assumed to be observed whenever the corresponding failure time is uncensored; when $\delta_i = 0$, $V_i$ is undefined and is not meaningful. The censoring time is assumed to be conditionally independent of $(T, V)$ given $Z$.

We consider a localized version of the log partial likelihood function for $\beta = \beta(v)$ at a fixed $v$:

$$l(v, \beta) = \sum_{i=1}^n \int_0^1 \int_0^\tau K_h(u - v) \left[\beta^T Z_i(t) - \log\left(\sum_{j=1}^n Y_j(t) e^{\beta^T Z_j(t)}\right)\right] \quad (3)$$
$$\times N_i(dt, du),$$

where $K_h(x) = K(x/h)/h$, $K(\cdot)$ is a kernel function with support $[-1, 1]$, $\tau$ is the end of the follow-up period and $h = h_n$ is a bandwidth. Here $Y_i(t) = I(X_i \geq t)$ and $N_i(t, v) = I(X_i \leq t, \delta_i = 1, V_i \leq v)$ is the marked point counting process with a jump at an uncensored failure times $X_i$ and the associated mark $V_i$. For background on marked point processes see Brémaud [2] and Martinussen and Scheike [10].

The log partial likelihood function (3) resembles that of Kalbfleisch and Prentice [8] in the case of discrete marks, except that it borrows strength from observations having marks in the neighborhood of $v$. The kernel function is designed to give greater weight to observations with marks near $v$ than those further away. The local maximum partial likelihood estimator of $\beta(v)$ is a maximizer $\hat{\beta}(v)$ of (3). A similar approach has been studied by Cai and Sun [3] for estimating time-dependent coefficients in Cox regression models.

Denote $\mu_j = \int u^j K(u)\,du$, $\nu_j = \int u^j K^2(u)\,du$ for $j = 0, 1, 2$. For $\beta \in \mathbb{R}^p$, $t \geq 0$, let

$$S^{(j)}(t,\beta) = n^{-1} \sum_{i=1}^{n} Y_i(t) \exp\{\beta^T Z_i(t)\} Z_i(t)^{\otimes j},$$

where for any $z \in \mathbb{R}^p$, we denote $z^{\otimes 0} = 1$, $z^{\otimes 1} = z$ and $z^{\otimes 2} = zz^T$. Define $s^{(j)}(t,\beta) = E S^{(j)}(t,\beta)$ and

$$J_n(t,\beta) = \frac{S^{(2)}(t,\beta)}{S^{(0)}(t,\beta)} - \left(\frac{S^{(1)}(t,\beta)}{S^{(0)}(t,\beta)}\right)^{\otimes 2},$$

$$J(t,\beta) = \frac{s^{(2)}(t,\beta)}{s^{(0)}(t,\beta)} - \left(\frac{s^{(1)}(t,\beta)}{s^{(0)}(t,\beta)}\right)^{\otimes 2}.$$

Taking the derivative of $l(v,\beta)$ with respect to $\beta$ gives the score function

(4)
$$\begin{aligned}U(v,\beta) &= l'_\beta(v,\beta) \\ &= \sum_{i=1}^{n} \int_0^1 \int_0^\tau K_h(u-v)\left[Z_i(t) - \frac{S^{(1)}(t,\beta)}{S^{(0)}(t,\beta)}\right] N_i(dt,du).\end{aligned}$$

The maximum partial likelihood estimator is a solution to $U(v,\hat{\beta}(v)) = 0$, and can be computed using a Newton–Raphson algorithm. The second derivative of $l(v,\beta)$ with respect to $\beta$ yields

$$l''_\beta(v,\beta) = -\sum_{i=1}^{n} \int_0^1 \int_0^\tau K_h(u-v) J_n(t,\beta) N_i(dt,du).$$

Although inference on $\beta$ is usually of primary interest, the baseline function $\lambda_0(t,v)$ can also be estimated, by smoothing the increments of the following estimator of the doubly cumulative baseline function $\Lambda_0(t,v) = \int_0^t \int_0^v \lambda_0(s,u)\,ds\,du$:

(5)
$$\hat{\Lambda}_0(t,v) = \int_0^t \int_0^v \frac{N(ds,du)}{nS^{(0)}(s,\hat{\beta}(u))}.$$

2.2. *Asymptotic results.* We make use of the following regularity conditions; not all of these conditions are required for the proof of each theorem, nor are they the minimum required set of conditions.

CONDITION A.

- (A.1) $\beta(v)$ has componentwise continuous second derivatives on $[0,1]$. The second partial derivative of $\lambda_0(t,v)$ with respect to $v$ exists and is continuous on $[0,\tau] \times [0,1]$. The covariate process $Z(t)$ has paths that are left-continuous and of bounded variation, and satisfies the moment condition $E[\|Z(t)\|^4 \exp(2M\|Z(t)\|)] < \infty$, where $M$ is a constant such that $(v, \beta(v)) \in [0,1] \times (-M, M)^p$ for all $v$ and $\|A\| = \max_{k,l} |a_{kl}|$ for a matrix $A = (a_{kl})$.
- (A.2) For $j = 0, 1, 2$, each component of $s^{(j)}(t, \theta)$ is continuous on $[0,\tau] \times [-M, M]^p$, and $\sup_{t \in [0,\tau], \theta \in [-M, M]^p} \|S^{(j)}(t,\theta) - s^{(j)}(t,\theta)\| = O_p(n^{-1/2})$.
- (A.3) $s^{(0)}(t, \theta) > 0$ on $[0,\tau] \times [-M, M]^p$ and the matrix $\Sigma(v) = \int_0^\tau J(t, \beta(v)) \times \lambda_0(t,v) s^{(0)}(t, \beta(v)) \, dt$ is positive definite.
- (A.4) $E(N_i(dt, dv)|\mathcal{F}_{t-}) = E(N_i(dt, dv)|Y_i(t), Z_i(t))$, where $\mathcal{F}_t = \sigma\{I(X_i \leq s, \delta_i = 1), I(X_i \leq s, \delta_i = 0), V_i I(X_i \leq s, \delta_i = 1), Z_i(s); 0 \leq s \leq t, i = 1, \ldots, n\}$ is the (right-continuous) filtration generated by $\{N_i(s,v), Y_i(s), Z_i(s); 0 \leq s \leq t, 0 \leq v \leq 1, i = 1, \ldots, n\}$.
- (A.5) The kernel function $K(\cdot)$ is symmetric with support $[-1, 1]$ and of bounded variation. The bandwidth satisfies $nh^2 \to \infty$ and $nh^5 \to 0$ as $n \to \infty$.

Note that the condition (A.2) holds under the condition (A.1) given some additional moment conditions on $Z(t) - Z(s)$ and $\exp(b^T Z(t)) - \exp(b^T Z(s))$. If $Z(t) = Z$, not depending on $t$, then (A.2) holds by the Donsker theorem (Theorem 19.5 of van der Vaart [19]). The condition (A.4) assumes that the mark-specific instantaneous failure rate at time $t$ given the observed information up to time $t$ only depends on the failure status and the current covariate value. Under (A.4) and by the definition (1), $E(N_i(dt, dv)|\mathcal{F}_{t-}) = Y_i(t)\lambda(t, v|Z_i(t)) \, dt \, dv$, and $M_i(t,v) = \int_0^t \int_0^v [N_i(ds, dx) - Y_i(s)\lambda(s, x|Z_i(s)) \, ds \, dx]$ is a martingale with respect to $\mathcal{F}_t$ for each fixed $v$ ([10], page 31). Further, it follows by Aalan and Johansen [1] that $M_i(\cdot, v_1)$ and $M_i(\cdot, v_2) - M_i(\cdot, v_1)$ are orthogonal square integrable martingales with respect to $\mathcal{F}_t$ for any $0 \leq v_1 \leq v_2 \leq 1$. To avoid the problems at the boundaries $v = 0, 1$, we shall study the asymptotic properties of $\hat{\beta}(v)$ for the interior values of $v \in [a, b] \subset (0, 1)$.

First we present the following result that is essential for proving the asymptotic normality of $\hat{\beta}(v)$ and provides insight into the constructions of the confidence bands and test statistics that follow. Let

$$(6) \quad \tilde{W}_A(v) = n^{-1/2} \sum_{i=1}^n \int_a^v \int_0^\tau A(u) \left[ Z_i(t) - \frac{s^{(1)}(t, \beta(u))}{s^{(0)}(t, \beta(u))} \right] M_i(dt, du),$$

where $A(u)$ is a deterministic $p \times p$ matrix with bounded components and $0 \leq a < b \leq 1$.

THEOREM 1. *Assume that each component of the $p \times p$ matrix $A(v)$, $v \in [a,b]$, is continuous. Under conditions (A.1)–(A.4), $\tilde{W}_A(v)$ converges weakly to a p-dimensional mean-zero Gaussian martingale, $W_A(v)$, with continuous sample paths on $v \in [a,b]$. The covariance matrix of $W_A(v)$ is given by $\operatorname{Cov}(W_A(v)) = \int_a^v A(u)\Sigma(u)A(u)\,du$.*

Let

$$(7) \qquad \hat{\Sigma}_{\hat{A}}(v) = n^{-1} \sum_{i=1}^n \int_a^v \int_0^\tau \hat{A}(u) J_n(t, \hat{\beta}(u)) \hat{A}^T(u) N_i(dt, du),$$

where $\hat{A}(v)$ is a consistent estimator of $A(v)$ uniformly in $v \in [a,b] \subset [0,1]$. It can be shown that $\hat{\Sigma}_{\hat{A}}(v)$ is a consistent estimator of $\operatorname{Cov}(W_A(v))$.

The consistency and asymptotic normality of $\hat{\beta}(v)$ are established in the next two theorems.

THEOREM 2. *Under conditions (A.1)–(A.5), $\hat{\beta}(v)$ converges to $\beta(v)$ uniformly in $v \in [a,b] \subset (0,1)$.*

THEOREM 3. *Under conditions (A.1)–(A.5), $(nh)^{1/2}(\hat{\beta}(v) - \beta(v)) \xrightarrow{\mathcal{D}} N(0, \nu_0 \Sigma^{-1}(v))$ for $v \in [a,b]$.*

The proof of Theorem 3 uses a Taylor expansion of the score function, leading to $\hat{\beta}(v) - \beta(v) = -(l''_\beta(v, \beta^*(v)))^{-1} U(\beta(v))$, where $\beta^*(v)$ is on the line segment between $\hat{\beta}(v)$ and $\beta(v)$. The asymptotic variance of $n^{-1/2} h^{1/2} U(\beta(v))$ is shown to be $\nu_0 \Sigma(v)$, which is the in probability limit of $\tilde{\Sigma}_n(\beta(v)) = n^{-1} h \times \sum_{i=1}^n \int_0^1 \int_0^\tau (K_h(u-v))^2 J_n(t, \beta(v)) N_i(dt, du)$. It can also be shown that $\hat{\Sigma}(v) \equiv -l''_\beta(v, \hat{\beta}(v))/n \xrightarrow{P} \Sigma(v)$ as $n \to \infty$. Thus, the asymptotic variance of $(nh)^{1/2} \times (\hat{\beta}(v) - \beta(v))$ can be estimated by $\hat{\Sigma}_1(v) = (l''_\beta(v, \hat{\beta}(v))/n)^{-1} \tilde{\Sigma}_n(\hat{\beta}(v)) (l''_\beta(v, \hat{\beta}(v))/n)^{-1}$. An alternative estimator is $\hat{\Sigma}_2(v) = -\nu_0 (l''_\beta(v, \hat{\beta}(v))/n)^{-1}$. It is easy to check that $\nu_0 = 3/5$ for Epanechnikov's kernel $K(x) = \frac{3}{4}(1-x^2)$, $-1 < x < 1$. Simulations indicate that the two estimators have similar finite sample performance.

Theorem 3 will lead to the construction of pointwise confidence intervals for $\operatorname{VE}(v)$. Simultaneous inference over $v \in [a,b]$ will be possible in terms of the estimate $\hat{B}(v) = \int_a^v \hat{\beta}(u)\,du$ of the cumulative regression coefficient $B(v) = \int_a^v \beta(u)\,du$. We have the following weak convergence result for $\hat{B}(v)$.

THEOREM 4. *Under conditions* (A.1)–(A.5), $n^{1/2}(\hat{B}(v) - B(v))$ *converges weakly to a p-dimensional mean-zero Gaussian martingale* $W_{\Sigma^{-1}}(v)$ *with continuous sample paths on* $v \in [a, b]$. *The covariance matrix of* $W_{\Sigma^{-1}}(v)$ *is* $\int_a^v \Sigma(u)^{-1} du$, *which can be consistently estimated by* $\hat{\Sigma}_{\hat{A}}(v)$ *defined by* (7) *with* $A(v) = (\Sigma(v))^{-1}$ *and* $\hat{A}(v) = (\hat{\Sigma}(v))^{-1}$.

2.3. *Confidence bands for vaccine efficacy.* Let $\beta(v) = (\beta_1(v), \beta_2^T(v))^T$. Then the vaccine efficacy can be expressed as $\mathrm{VE}(v) = 1 - \exp(\beta_1(v))$. The estimated vaccine efficacy is $\widehat{\mathrm{VE}}(v) = 1 - \exp(\hat{\beta}_1(v))$. By Theorem 3 and the delta method, $(nh)^{1/2}(\widehat{\mathrm{VE}}(v) - \mathrm{VE}(v)) \xrightarrow{\mathcal{D}} N(0, \nu_0 \sigma_1^2(v) \exp(2\beta_1(v)))$ for $v \in [a, b]$, where $\sigma_1^2(v)$ is the first element on the diagonal of $\Sigma^{-1}(v)$. Let $\hat{\sigma}^2_{\beta_1}(v)$ be the first element on the diagonal of $\hat{\Sigma}_1(v)$. By the discussions on the consistent estimators for the asymptotic variance following Theorem 3, $\hat{\sigma}^2_{\beta_1}(v)$ is a consistent estimator for $\nu_0 \sigma_1^2(v)$. A pointwise $100(1-\alpha)\%$ confidence band for $\mathrm{VE}(v)$ is given by

$$(8) \qquad \widehat{\mathrm{VE}}(v) \pm (nh)^{-1/2} z_{\alpha/2} \hat{\sigma}_{\beta_1}(v) \exp(\hat{\beta}_1(v)), \qquad a \leq v \leq b,$$

where $z_{\alpha/2}$ is the upper $\alpha/2$ quantile of the standard normal distribution.

To derive simultaneous confidence bands for the cumulative vaccine efficacy $\mathrm{CV}(v) = \int_a^v \mathrm{VE}(u)\, du$, we consider the point estimator $\widehat{\mathrm{CV}}(v) = \int_a^v \widehat{\mathrm{VE}}(u)\, du$. Then

$$\sqrt{n}(\widehat{\mathrm{CV}}(v) - \mathrm{CV}(v)) = \sqrt{n} \int_a^v (\exp(\beta_1(v)) - \exp(\hat{\beta}_1(v)))\, du.$$

Note that $\sqrt{n}(\widehat{\mathrm{CV}}(v) - \mathrm{CV}(v)) \approx \sqrt{n} \int_a^v \exp(\beta_1(v)(\beta_1(v) - \hat{\beta}_1(v))\, du$. From the proof of Theorem 4, it can be shown that $\sqrt{n}(\widehat{\mathrm{CV}}(v) - \mathrm{CV}(v))$ converges weakly to a mean-zero Gaussian process, $e_1^T W_A(v)$, $a \leq v \leq b$, with continuous paths and independent increments, where $A(v) = \exp(\beta_1(v))\Sigma(v)^{-1}$ and $e_1$ is the first column of the $p \times p$ identity matrix. The variance of $e_1^T W_A(v)$ equals $\rho^2(v) = \int_a^v \sigma_1^2(u) \exp(2\beta_1(u))\, du$ by Theorem 1, which can be conveniently estimated by $\int_a^v \hat{\sigma}_1^2(u) \exp(2\hat{\beta}_1(u))\, du$, where $\hat{\sigma}_1^2(v)$ is the first element of the diagonal of $\hat{\Sigma}(v)^{-1}$. We suspect that this estimator may ignore the finite sample correlations of $\beta_1(v) - \hat{\beta}_1(v)$ at different values of $v$, thus over- or underestimating the true variance. We propose to use $\hat{\rho}^2(v) = e_1^T \hat{\Sigma}_{\hat{A}}(v) e_1$ as the estimator of the asymptotic variance of $\sqrt{n}(\widehat{\mathrm{CV}}(v) - \mathrm{CV}(v))$, where $\hat{\Sigma}_{\hat{A}}(v)$ is obtained from (7) with $\hat{A}(v) = \exp(\hat{\beta}_1(v))\hat{\Sigma}(v)^{-1}$, which is uniformly consistent by Theorem 1. Consequently, a pointwise $100(1-\alpha)\%$ confidence band for $\mathrm{CV}(v)$ is given by

$$(9) \qquad \widehat{\mathrm{CV}}(v) \pm n^{-1/2} z_{\alpha/2} \hat{\rho}(v), \qquad a \leq v \leq b.$$

Let $\mathcal{V}$ be a set of values of $v$ in $[a,b]$. We may take $\mathcal{V} = [a,b]$ or $\mathcal{V} = \{v_k, k = 1,\ldots,K\}$ with $v_1 < \cdots < v_K$. Note that if $U(v)$ is a Gaussian martingale with variance $\rho^2(v)$, for $a \leq v \leq b$, then $U(v)\rho(b)[\rho^2(b) + \rho^2(v)]^{-1}$ has the same distribution as $B^0(\rho^2(v)/(\rho^2(b) + \rho^2(v)))$, $a \leq v \leq b$, where $B^0(\cdot)$ is a Brownian bridge. By the weak convergence of $\sqrt{n}(\widehat{\mathrm{CV}}(v) - \mathrm{CV}(v))$, the uniform consistency of $\hat{\rho}^2(v)$ to $\rho^2(v)$ and the continuous mapping theorem, we have

$$\sup_{v \in \mathcal{V}}|\sqrt{n}(\widehat{\mathrm{CV}}(v) - \mathrm{CV}(v))\hat{\rho}(b)/(\hat{\rho}^2(b) + \hat{\rho}^2(v))|$$

$$\xrightarrow{\mathcal{D}} \sup_{v \in \mathcal{V}}|B^0(\rho^2(v)/(\rho^2(b) + \rho^2(v)))|.$$

Thus a simultaneous $100(1-\alpha)\%$ confidence band for $\mathrm{CV}(v)$, $v \in \mathcal{V}$, is given by

(10) $$\widehat{\mathrm{CV}}(v) \pm n^{-1/2} u_\alpha [\hat{\rho}^2(b) + \hat{\rho}^2(v)]/\hat{\rho}(b),$$

where $u_\alpha$ is the upper $\alpha$-quantile of the distribution of $\sup_{v \in \mathcal{V}} |B^0(\rho^2(v)/(\rho^2(b) + \rho^2(v)))|$. The $u_\alpha$ is the upper $\alpha$-quantile of $\sup_{0 \leq v \leq 1/2} |B^0(v)|$ if $\mathcal{V} = [a,b]$, which has been tabulated by Schumacher [15] for some $\alpha$ values. In the simulation study presented in the next section, we estimate $u_\alpha$ by the upper $\alpha$-quantile of the distribution of $\sup_{v_k \in \mathcal{V}} |B^0(\hat{\rho}^2(v_k)/(\hat{\rho}^2(b) + \hat{\rho}^2(v_k)))|$ in both cases when $\mathcal{V} = [a,b]$ or $\mathcal{V} = \{v_k, k = 1,\ldots,K\}$, which can be obtained by simulating a Brownian bridge for given $\hat{\rho}^2(v)$.

Alternatively, other resampling techniques such as the Gaussian multiplier method of Lin, Wei and Ying [9] can be used to estimate the critical value $u_\alpha$. This method can be briefly outlined as follows. Let $\xi_1,\ldots,\xi_n$ be i.i.d. standard normal random variables and

(11) $$W^*_{\hat{A}}(v) = n^{-1/2} \sum_{i=1}^n \xi_i \int_0^v \int_0^\tau \hat{A}(u)\left[Z_i(t) - \frac{S^{(1)}(t,\hat{\beta}(u))}{S^{(0)}(t,\hat{\beta}(u))}\right] M_i(dt,du).$$

Then the distribution $\sqrt{n}(\widehat{\mathrm{CV}}(v) - \mathrm{CV}(v))$ can be approximated by the conditional distribution of $e_1^T W^*_{\hat{A}}(v)$ given the observed data sequence, where $\hat{A} = \exp(\hat{\beta}_1(v)) \times (\hat{\Sigma}(v))^{-1}$. Consequently, the distribution of $\sup_{v \in \mathcal{V}} |\sqrt{n}(\widehat{\mathrm{CV}}(v) - \mathrm{CV}(v))\hat{\rho}(b)[\hat{\rho}^2(b) + \hat{\rho}^2(v)]^{-1}|$ can be approximated by the conditional distribution of $U^* = \sup_{v \in \mathcal{V}} |e_1^T W^*_{\hat{A}}(v)\hat{\rho}(b)[\hat{\rho}^2(b) + \hat{\rho}^2(v)]^{-1}|$ given the observed data sequence. Let $u^*_\alpha$ be the $(1-\alpha)$-quantile of the copies of $U^*$ obtained by repeatedly generating sets of i.i.d. standard normal random variables. A simultaneous $100(1-\alpha)\%$ confidence band for $\mathrm{CV}(v)$, $v \in \mathcal{V}$, is given by

(12) $$\widehat{\mathrm{CV}}(v) \pm n^{-1/2} u^*_\alpha [\hat{\rho}^2(b) + \hat{\rho}^2(v)]/\hat{\rho}(b).$$

This resampling technique is also applicable to the hypothesis tests for vaccine efficacy developed in the next subsection.

2.4. *Testing vaccine efficacy.* We are interested in testing the following two sets of hypotheses. The first set of hypotheses is

$$H_{10}: \text{VE}(v) = 0 \text{ for } v \in [a, b]$$

versus $\quad H_{1a}: \text{VE}(v) \neq 0$ for some $v$ (general alternative)

or $\quad H_{1m}: \text{VE}(v) \geq 0$ with strict inequality for at least some $v$

(monotone alternative).

The second set of hypotheses is

$$H_{20}: \text{VE}(v) \text{ does not depend on } v \in [a, b]$$

versus $\quad H_{2a}: \text{VE}(v)$ depends on $v$ (general alternative)

or $\quad H_{2m}: \text{VE}(v)$ decreases as $v$ increases (monotone alternative).

Let $\beta_1(v)$ be the first component of $\beta(v)$. Then the null hypothesis $H_{10}$ is equivalent to $\beta_1(v) = 0$ and the null hypothesis $H_{20}$ is equivalent to $\beta_1(v)$ does not depend on $v$. The null hypothesis $H_{10}$ implies the vaccine affords no protection against any infecting strain of virus. The alternative $H_{1m}$ indicates that the vaccine provides protection for at least some of the infecting strains, while $H_{1a}$ states that the vaccine provides either protection or increased risk for some infecting strains. The null hypothesis $H_{20}$ implies there is no difference in vaccine effect for different infecting strains, measured by their distance $v$ to the strains contained in the vaccine. The ordered alternative $H_{2m}$ states that vaccine efficacy decreases with $v$ and the alternative $H_{2a}$ indicates that the vaccine efficacy changes with $v$.

In this section, we develop some test procedures for detecting departures from $H_{10}$ in the direction of $H_{1m}$ and $H_{1a}$ and for detecting departures from $H_{20}$ in the direction of $H_{2m}$ and $H_{2a}$. By Theorem 4 and the discussions in Section 2.3, the process $\sqrt{n}(\widehat{\text{CV}}(v) - \text{CV}(v))$, $a \leq v \leq b$, converges weakly to a Gaussian martingale with predictable variation $\rho^2(v)$. Let $\xi(v) = \sqrt{n}(\widehat{\text{CV}}(v) - \text{CV}(v))/\rho(b)$. It follows from Theorem 4 that $\xi(v) \xrightarrow{\mathcal{D}} W(t(v))$, $a \leq v \leq b$, where $W(\cdot)$ is a Wiener process and $t(v) = \rho^2(v)/\rho^2(b)$.

To test $H_{10}$, let $\hat{Z}^{(1)}(v) = \sqrt{n}\widehat{\text{CV}}(v)/\hat{\rho}(b)$ and $\hat{t}(v) = \hat{\rho}^2(v)/\hat{\rho}^2(b)$. Consider the following test statistics:

$$T_a^{(1)} = \int_a^b (\hat{Z}^{(1)}(v))^2 \, d\hat{t}(v), \qquad T_{m1}^{(1)} = \int_a^b \hat{Z}^{(1)}(v) \, d\hat{t}(v).$$

These test statistics have somewhat complicated null distributions (see below) so we also consider the following simpler test statistic based on a finite grid, which leads to a standard normal null distribution:

$$T_{m2}^{(1)} = (K-1)^{-1/2} \sum_{k=2}^{K} (\hat{Z}^{(1)}(v_k) - \hat{Z}^{(1)}(v_{k-1}))/(\hat{t}(v_k) - \hat{t}(v_{k-1}))^{1/2},$$

where $a \leq v_1 < \cdots < v_K \leq b$ are the grid points in $[a,b]$. A similar test statistic with a standard normal null distribution is also proposed for $H_{20}$ later. Under $H_{10}$, $T_a^{(1)} \xrightarrow{\mathcal{D}} \int_a^b (W(t(v)))^2 \, dt(v) \stackrel{\mathcal{D}}{=} \int_0^1 (W(t))^2 \, dt$, $T_{m1}^{(1)} \xrightarrow{\mathcal{D}} \int_a^b W(t(v)) \, dt(v) \stackrel{\mathcal{D}}{=} \int_0^1 W(t) \, dt$ and $T_{m2}^{(1)} \xrightarrow{\mathcal{D}} N(0,1)$. The distributions of $T_a^{(1)}$ and $T_{m1}^{(1)}$ under $H_{10}$ can also be approximated by those of $\int_a^b (W(\hat{t}(v)))^2 \, d\hat{t}(v)$ and $\int_a^b W(\hat{t}(v)) \, d\hat{t}(v)$ for given $\hat{t}(v)$, respectively, which are used in the numerical studies for better finite sample approximations. We denote the upper $\alpha$-quantiles of these two distributions by $c_a^{(1)}$ and $c_{m1}^{(1)}$, respectively.

The test statistic $T_a^{(1)}$ captures general departures $H_{1a}$, while the test statistics $T_{m1}^{(1)}$ and $T_{m2}^{(1)}$ are sensitive to the monotone departure $H_{1m}$. Both test statistics $T_{m1}^{(1)}$ and $T_{m2}^{(1)}$ are likely to be positive when $\text{VE}(v) \geq 0$ for all $v$ with strict inequality for some $v$. Hence the tests based on $T_a^{(1)}$, $T_{m1}^{(1)}$ and $T_{m2}^{(1)}$ reject $H_{10}$ if $T_a^{(1)} > c_a^{(1)}$, $T_{m1}^{(1)} > c_{m1}^{(1)}$ and $T_{m2}^{(1)} > z_\alpha$, respectively.

To test $H_{20}$, let $\hat{Z}^{(2)}(v) = \sqrt{n}(\frac{1}{v-a}\widehat{\text{CV}}(v) - \frac{1}{b-a}\widehat{\text{CV}}(b))/\hat{\rho}(b)$. Note that, under $H_{20}$, $\hat{Z}^{(2)}(v) = \sqrt{n}[\frac{1}{v-a}(\widehat{\text{CV}}(v) - \text{CV}(v)) - \frac{1}{b-a}(\widehat{\text{CV}}(b) - \text{CV}(b))]/\hat{\rho}(b)$. By Theorem 4 and the continuous mapping theorem, under $H_{20}$, $\hat{Z}^{(2)}(v) \xrightarrow{\mathcal{D}} \frac{1}{v-a}W(t(v)) - \frac{1}{b-a}W(1) \equiv Z^{(2)}(v)$ for $v \in [a_1, b]$, where $a < a_1 < b$. We propose the following test statistics for evaluating $H_{20}$:

$$T_a^{(2)} = \int_{a_1}^b (\hat{Z}^{(2)}(v))^2 \, d\hat{t}(v), \qquad T_{m1}^{(2)} = \int_{a_1}^b \hat{Z}^{(2)}(v) \, d\hat{t}(v),$$

$$T_{m2}^{(2)} = \hat{\Pi}_K^{-1} \sum_{k=2}^K (\hat{Z}^{(2)}(v_{k-1}) - \hat{Z}^{(2)}(v_k))/\hat{\pi}_k,$$

where $a_1 \leq v_1 < \cdots < v_K \leq b$ are $K$ grid points in $[a_1, b]$, $\hat{\pi}_k^2$ is an estimate of the variance $\pi_k^2 = \text{Var}(Z^{(2)}(v_{k-1}) - Z^{(2)}(v_k))$ and $\hat{\Pi}_K^2$ is an estimate of the variance $\Pi_K^2$ of $\sum_{k=2}^K (Z^{(2)}(v_{k-1}) - Z^{(2)}(v_k))/\pi_k$. By the covariance of the Wiener process, it is easy to show that

$$\tau_{i,j} = \text{Cov}(Z^{(2)}(v_i), Z^{(2)}(v_j))$$
$$= \frac{t(v_i)}{(v_i - a)(v_j - a)} - \frac{t(v_i)}{(v_i - a)(b - a)} - \frac{t(v_j)}{(v_j - a)(b - a)} + \frac{1}{(b - a)^2},$$

for $v_i \leq v_j$. Thus, $\pi_k^2 = \tau_{k-1,k-1} - 2\tau_{k-1,k} + \tau_{k,k}$. Let $\Gamma = (\tau_{i,j})_{K \times K}$ and

$$\xi^T = (\pi_2^{-1}, \pi_3^{-1} - \pi_2^{-1}, \ldots, \pi_K^{-1} - \pi_{K-1}^{-1}, -\pi_K^{-1}).$$

It follows that $\Pi_K = \xi^T \Gamma \xi$. The estimates $\hat{\pi}_k^2$ and $\hat{\Pi}_K^2$ are obtained by replacing $t(v)$ with $\hat{t}(v)$.

By the weak convergence of $\hat{Z}^{(2)}(v)$ to $Z^{(2)}(v)$, and the convergence in probability of $\hat{t}(v)$ to $t(v)$, $a_1 \leq v \leq b$, we have $T_{m2}^{(2)} \xrightarrow{\mathcal{D}} N(0,1)$ under $H_{20}$. It also follows that $T_a^{(2)} \xrightarrow{\mathcal{D}} \int_{a_1}^{b} (Z^{(2)}(v))^2 \, dt(v)$, and $T_{m1}^{(2)} \xrightarrow{\mathcal{D}} \int_{a_1}^{b} Z^{(2)}(v) \, dt(v)$ under $H_{20}$. The distributions of $T_a^{(2)}$ and $T_{m1}^{(2)}$ under $H_{20}$ can be approximated by those of $\int_{a_1}^{b} (W(\hat{t}(v))/(v-a) - W(\hat{t}(b))/(b-a))^2 \, d\hat{t}(v)$ and $\int_{a_1}^{b} (W(\hat{t}(v))/(v-a) - W(\hat{t}(b))/(b-a)) \, d\hat{t}(v)$ for given $\hat{t}(v)$, respectively, which are used in the numerical studies for better finite sample approximations. We denote the upper $\alpha$-quantiles of these two distributions by $c_a^{(2)}$ and $c_{m1}^{(2)}$, respectively.

The test statistic $T_a^{(2)}$ captures general departures $H_{2a}$ while the test statistics $T_{m1}^{(2)}$ and $T_{m2}^{(2)}$ are sensitive to the monotone departure $H_{2m}$. Both $T_{m1}^{(2)}$ and $T_{m2}^{(2)}$ are expected to be positive when $\text{VE}(v)$ decreases as $v$ increases, that is, when $H_{2m}$ holds. Hence the tests $T_a^{(2)}$, $T_{m1}^{(2)}$ and $T_{m2}^{(2)}$ reject $H_{20}$ if $T_a^{(2)} > c_a^{(2)}$, $T_{m1}^{(2)} > c_{m1}^{(2)}$ and $T_{m2}^{(2)} > z_\alpha$, respectively.

**3. Simulation study.** In this section, we conduct a simulation study to check the finite sample performance of the proposed estimation and hypothesis testing procedures using the simple mark-specific proportional hazards model:

(13) $$\lambda(t, v|z) = \exp\{\gamma v + (\alpha + \beta v)z\}, \qquad t \geq 0, 0 \leq v \leq 1,$$

where $\alpha$, $\beta$ and $\gamma$ are constants and the treatment indicator $z$ takes value 0 or 1 with probability of 0.5 for each value. Under model (13), the mark-specific baseline function is $\lambda_0(t, v) = \exp(\gamma v)$ and $\text{VE}(v) = 1 - \exp(\alpha + \beta v)$. The null hypothesis $H_{10}$ of no vaccine efficacy holds if both $\alpha = 0$ and $\beta = 0$, and the null hypothesis $H_{20}$ that vaccine efficacy does not depend on the type of infecting strain is true if $\beta = 0$. Various choices of $\alpha$ and $\beta$ specify different alternatives for $H_{10}$ and $H_{20}$.

We consider the following simulation models:

(M1) $(\alpha, \beta, \gamma) = (0, 0, 0.3)$, for the null hypothesis $H_{10}$ of no vaccine efficacy;
(M2) $(\alpha, \beta, \gamma) = (-0.5, 0.5, 0.3)$, as the first alternative of $H_{10}$;
(M3) $(\alpha, \beta, \gamma) = (-0.6, 0.6, 0.3)$, as the second alternative of $H_{10}$;
(M4) $(\alpha, \beta, \gamma) = (-0.6, 0, 0.3)$, as the third alternative of $H_{10}$;
(M5) $(\alpha, \beta, \gamma) = (-0.69, 0, 0.3)$, for the null hypothesis $H_{20}$ that vaccine efficacy does not depend on the type of infecting strain;
(M6) $(\alpha, \beta, \gamma) = (-1.2, 1.2, 0.3)$, as the first alternative of $H_{20}$;
(M7) $(\alpha, \beta, \gamma) = (-1.5, 1.5, 0.3)$, as the second alternative of $H_{20}$;
(M8) $(\alpha, \beta, \gamma) = (-1.8, 1.8, 0.3)$, as the third alternative of $H_{20}$.

The models (M2) to (M4) are considered as the alternatives for $H_{1m}$ and $H_{1a}$. The departure from $H_{10}: \text{VE}(v) = 0$ increases as the simulation model moves from (M2) to (M4). The models (M6) to (M8) are considered as the alternatives for $H_{2m}$ and $H_{2a}$. The departure from $H_{20}$ increases as the simulation model moves from (M6) to (M8).

We generate the censoring times from an exponential distribution, independent of $(T, V)$, with the censoring rates ranging from 20% to 30%. We set the interval of analyses for $v$ as $[a, b] = [0.1, 0.9]$ and bandwidths are chosen as $h = 0.05, 0.1, 0.15$. The observed failure times with marks outside the interval $[a, b]$ can also be used since the smoothing at $v$ takes the cases with marks in its $h$-neighborhood. The Epanechnikov kernel $K(x) = 0.75(1 - x^2)I\{|x| \leq 1\}$ is used throughout. Sample sizes of $n = 500$ and 800 are studied.

For the tests $T_{m2}^{(1)}$ and $T_{m2}^{(2)}$, we take the grid of eight evenly spaced points in $[a, b]$ from 0.196 to 0.868. Table 1 lists the empirical sizes and powers of the test statistics $T_a^{(1)}$, $T_{m1}^{(1)}$ and $T_{m2}^{(1)}$ and Table 2 for the test statistics $T_a^{(2)}$, $T_{m1}^{(2)}$ and $T_{m2}^{(2)}$. The significance levels of these tests are given at $\alpha = 0.05$. Both tables also list the coverage probabilities of the 95% simultaneous confidence intervals for CV$(v)$, for $v \in [a, b]$ and for $v$ in the grid. The critical values for the tests $T_{m2}^{(1)}$ and $T_{m2}^{(2)}$ at nominal level 0.05 are $z_\alpha = 1.645$. The critical values for the tests $T_a^{(1)}$ and $T_a^{(2)}$, $T_{m1}^{(1)}$ and $T_{m1}^{(2)}$ are obtained by generating 10,000 Wiener processes $W(\cdot)$ with time parameter equal to $\hat{t}(v)$ and calculating the corresponding functionals of $W(\hat{t}(v))$, as described in the previous section. Each entry in Tables 1 and 2 is based on 1000 repetitions.

Most tests have appropriate sizes close to 5%. The test $T_a^{(2)}$ seems to be conservative for the simulation models used in the study. The test $T_{m1}^{(1)}$ has better power than the tests $T_a^{(1)}$ and $T_{m2}^{(1)}$. The test $T_{m1}^{(2)}$ has better power than the tests $T_a^{(2)}$ and $T_{m2}^{(2)}$. Therefore the tests that incorporate $\widehat{\text{CV}}(v)$ over the entire range $[a, b]$ present greater power than the simpler tests based on $\widehat{\text{CV}}(v)$ over the grid. We also observed that the powers of the tests seem to be influenced by the selection of bandwidth, with greater power for a larger bandwidth. Similar plots (not included here) to Figure 1 and Figure 2 but with larger bandwidth $h = 0.2$ show that the estimated standard errors of $\widehat{\text{CV}}(v)$ become smaller for larger $h$ while the biases stay approximately the same, resulting in increased power for the larger bandwidth. We suspect that this phenomenon is associated with the sample size and the convergence rate of the normalized $\widehat{\text{CV}}(v)$ to a Wiener process. The dependence of the power on the bandwidth should become small as the sample size increases. Further study on the bandwidth selection is warranted.

The coverage probabilities of the simultaneous confidence intervals for CV$(v)$ are closer to the 95% nominal level for $v$ on the grid than on $[a, b]$. This may be explained by the fact that the convergence for $v$ over the entire range $[a, b]$ is slower than the convergence on the grid. The evaluations of the proposed estimators for $\beta(v)$, VE$(v)$ and CV$(v)$ and their respective estimators of the standard deviations under some of the simulation models are presented in Figure 1 and Figure 2. The plots of the pointwise coverage probabilities for $\mathrm{VE}(v)$ and for $\mathrm{CV}(v)$ are given in Figure 3. These plots are based on $n=500$ and $h=0.1$.

TABLE 1
*Empirical sizes and powers of the tests $T_a^{(1)}$, $T_{m1}^{(1)}$ and $T_{m2}^{(1)}$ at the nominal level 0.05, and coverage probabilities of the 95% simultaneous confidence intervals for $\mathrm{CV}(v)$ with $v$ on the grid and on $[a,b]$*

| Model | $(\alpha,\beta,\gamma)$ | $n$ | $h$ | Size/Power | | | Coverage | |
|---|---|---|---|---|---|---|---|---|
| | | | | $T_a^{(1)}$ | $T_{m1}^{(1)}$ | $T_{m2}^{(1)}$ | Grid | $[a,b]$ |
| M1 | $(0,0,0.3)$ | 500 | 0.05 | 2.9 | 3.1 | 7.8 | 97.5 | 98.1 |
| | | | 0.1 | 4.9 | 5.9 | 8.3 | 96.6 | 97.4 |
| | | | 0.15 | 5.1 | 6.9 | 7.3 | 96.2 | 96.8 |
| | | 800 | 0.05 | 5.3 | 2.8 | 6.9 | 95.9 | 96.8 |
| | | | 0.1 | 5.7 | 4.7 | 6.8 | 95.5 | 97.0 |
| | | | 0.15 | 5.8 | 5.2 | 6.3 | 95.6 | 96.5 |
| M2 | $(-0.5, 0.5, 0.3)$ | 500 | 0.05 | 45.4 | 56.3 | 63.2 | 97.6 | 98.0 |
| | | | 0.1 | 60.3 | 71.4 | 65.7 | 97.0 | 97.5 |
| | | | 0.15 | 66.0 | 77.4 | 65.5 | 96.7 | 97.6 |
| | | 800 | 0.05 | 69.1 | 78.4 | 77.5 | 96.1 | 96.8 |
| | | | 0.1 | 80.3 | 86.5 | 80.1 | 95.6 | 96.7 |
| | | | 0.15 | 82.9 | 89.1 | 80.1 | 96.0 | 97.2 |
| M3 | $(-0.6, 0.6, 0.3)$ | 500 | 0.05 | 59.7 | 70.0 | 76.5 | 97.5 | 98.0 |
| | | | 0.1 | 75.4 | 83.9 | 78.8 | 96.9 | 97.8 |
| | | | 0.15 | 80.9 | 87.2 | 78.5 | 96.9 | 97.9 |
| | | 800 | 0.05 | 83.7 | 90.4 | 87.6 | 96.2 | 96.9 |
| | | | 0.1 | 90.8 | 94.4 | 89.6 | 96.0 | 96.8 |
| | | | 0.15 | 93.0 | 96.0 | 89.6 | 96.2 | 97.2 |
| M4 | $(-0.6, 0, 0.3)$ | 500 | 0.05 | 96.0 | 95.6 | 99.9 | 97.0 | 97.8 |
| | | | 0.1 | 99.1 | 98.8 | 100 | 96.7 | 97.6 |
| | | | 0.15 | 99.5 | 99.7 | 100 | 96.7 | 97.4 |
| | | 800 | 0.05 | 99.9 | 99.5 | 100 | 97.0 | 98.0 |
| | | | 0.1 | 100 | 100 | 100 | 96.9 | 97.3 |
| | | | 0.15 | 100 | 100 | 100 | 96.4 | 97.4 |

Now we demonstrate with a simulation example that the adoption of a standard method for testing the vaccine efficacy that ignores the mark is inefficient and can be misleading. We consider a special case of the model discussed in the Introduction, with $\lambda(t,v|z=0)=1$ and $\lambda(t,v|z=1)=2v$, for $t \geq 0$, $0 \leq v \leq 1$. The covariate $z$ is again a treatment indicator taking values 0 and 1 with probability of 0.5 for each value. The marginal hazards model ignoring the mark is therefore $\lambda(t|z=0)=1$ and $\lambda(t|z=1)=1$, for $t \geq 0$. The rest of the simulation setup such as the percentage of censorship, the kernel function and the bandwidth is the same as for the previous models. The model considered here represents both a proportional mark-

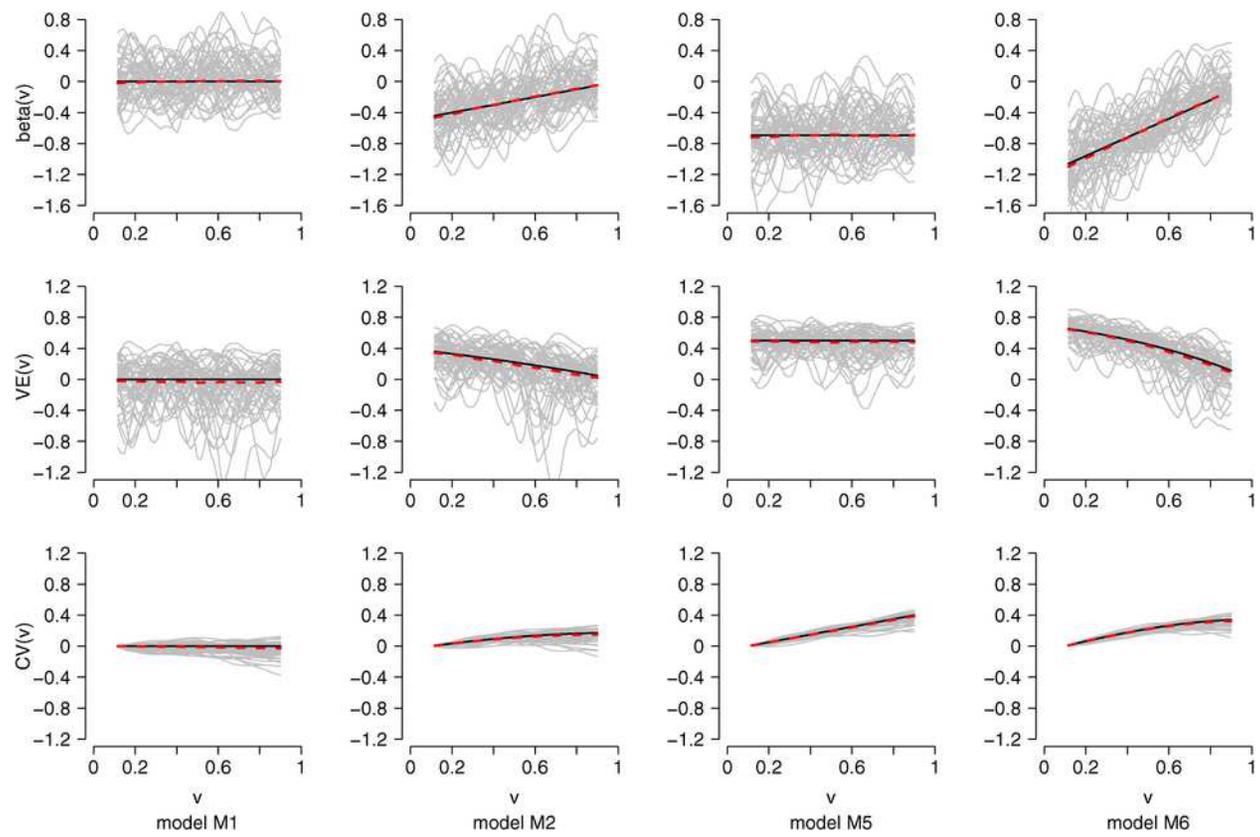

FIG. 1. *Plots of estimates for $\beta(v)$, $\text{VE}(v)$ and $\text{CV}(v)$ under the models M1, M2, M5 and M6 for $n = 500$, $h = 0.1$. The solid dark lines are the true functions and the dashed lines are the averages of the estimates based on 1000 repetitions. The gray lines are the corresponding estimates for $\beta(v)$, $\text{VE}(v)$ and $\text{CV}(v)$ of 50 random samples.*

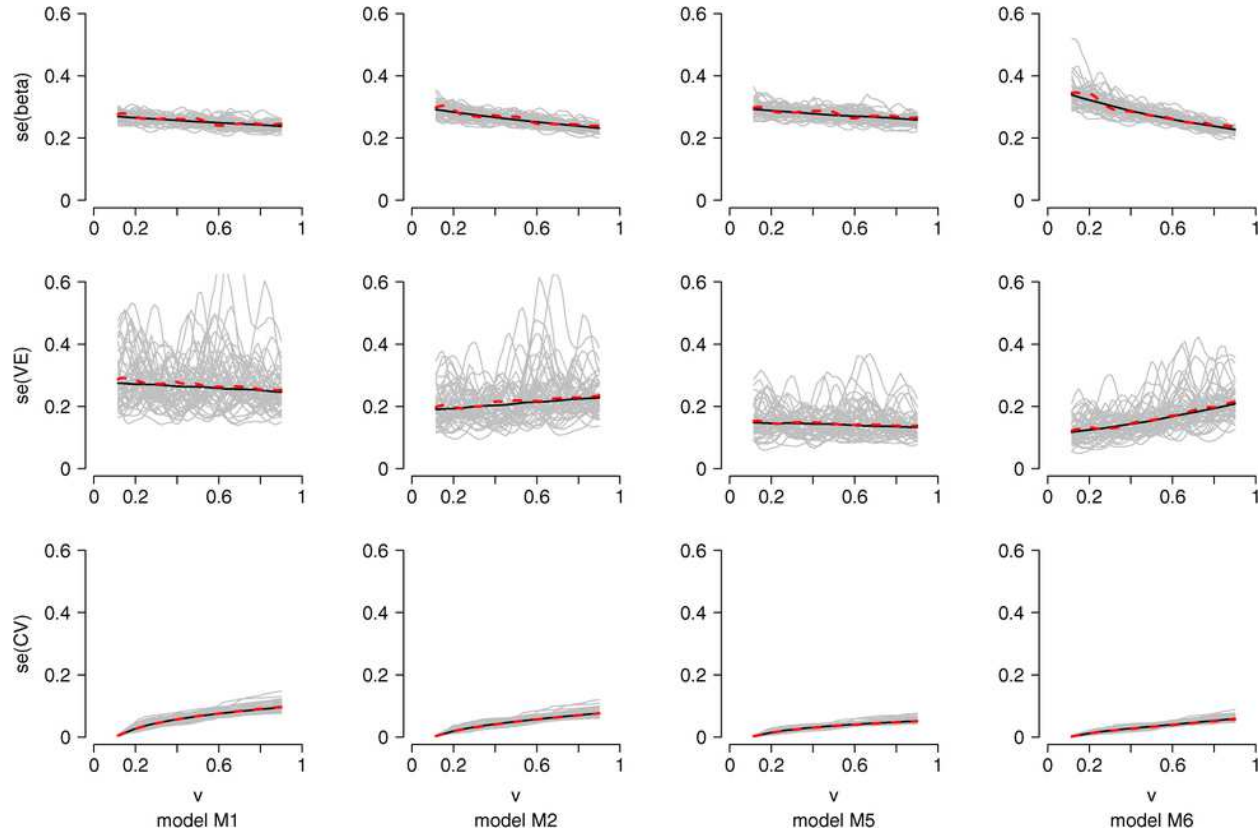

FIG. 2. *Plots of the standard errors under the models M1, M2, M5 and M6, based on $n = 500$, $h = 0.1$. The solid lines are the averages of the estimates of the standard deviations of $\hat{\beta}(v)$, $\widehat{\text{VE}}(v)$ and $\widehat{\text{CV}}(v)$, while the dashed lines are the sample standard deviations of $\hat{\beta}(v)$, $\widehat{\text{VE}}(v)$ and $\widehat{\text{CV}}(v)$, based on 1000 repetitions. The gray lines are the corresponding estimates for the standard deviations of $\hat{\beta}(v)$, $\widehat{\text{VE}}(v)$ and $\widehat{\text{CV}}(v)$ of 50 random samples.*

MARK-SPECIFIC PROPORTIONAL HAZARDS MODELS 17TABLE 2
*Empirical sizes and powers of the tests $T_a^{(2)}$, $T_{m1}^{(2)}$ and $T_{m2}^{(2)}$ at the nominal level 0.05, and coverage probabilities of the 95% simultaneous confidence intervals for $\mathrm{CV}(v)$ with $v$ on the grid and on $[a,b]$*

| | | | | Size/Power | | | Coverage | |
|---|---|---|---|---|---|---|---|---|
| **Model** | $(\alpha,\beta,\gamma)$ | $n$ | $h$ | $T_a^{(2)}$ | $T_{m1}^{(2)}$ | $T_{m2}^{(2)}$ | Grid | $[a,b]$ |
| M5 | $(-0.69, 0, 0.3)$ | 500 | 0.05 | 1.6 | 3.7 | 3.7 | 97.0 | 97.8 |
| | | | 0.1 | 2.1 | 3.7 | 4.5 | 96.5 | 97.5 |
| | | | 0.15 | 2.1 | 3.5 | 4.6 | 96.8 | 97.3 |
| | | 800 | 0.05 | 2.3 | 4.0 | 2.9 | 97.3 | 98.3 |
| | | | 0.1 | 2.6 | 4.3 | 3.2 | 97.0 | 97.6 |
| | | | 0.15 | 2.1 | 3.5 | 3.0 | 96.9 | 97.4 |
| M6 | $(-1.2, 1.2, 0.3)$ | 500 | 0.05 | 47.2 | 67.6 | 47.7 | 97.9 | 98.5 |
| | | | 0.1 | 60.2 | 76.7 | 62.3 | 97.1 | 97.6 |
| | | | 0.15 | 63.2 | 80.3 | 73.3 | 97.5 | 97.8 |
| | | 800 | 0.05 | 69.2 | 85.1 | 69.2 | 96.5 | 97.2 |
| | | | 0.1 | 80.4 | 92.0 | 80.4 | 96.6 | 97.6 |
| | | | 0.15 | 84.2 | 94.1 | 88.4 | 96.9 | 97.8 |
| M7 | $(-1.5, 1.5, 0.3)$ | 500 | 0.05 | 63.8 | 81.4 | 62.1 | 97.7 | 98.0 |
| | | | 0.1 | 76.9 | 78.0 | 63.6 | 97.2 | 98.0 |
| | | | 0.15 | 81.2 | 91.7 | 86.3 | 97.6 | 98.0 |
| | | 800 | 0.05 | 85.1 | 94.4 | 82.6 | 96.2 | 97.1 |
| | | | 0.1 | 93.2 | 98.2 | 91.8 | 96.1 | 97.6 |
| | | | 0.15 | 96.0 | 98.9 | 97.4 | 96.7 | 97.7 |
| M8 | $(-1.8, 1.8, 0.3)$ | 500 | 0.05 | 77.6 | 89.1 | 73.6 | 97.8 | 98.5 |
| | | | 0.1 | 87.1 | 95.6 | 85.7 | 97.3 | 98.4 |
| | | | 0.15 | 91.5 | 96.9 | 92.8 | 97.7 | 98.7 |
| | | 800 | 0.05 | 93.5 | 98.2 | 91.4 | 96.4 | 97.4 |
| | | | 0.1 | 98.2 | 99.5 | 97.0 | 96.3 | 97.5 |
| | | | 0.15 | 99.3 | 99.9 | 99.2 | 96.5 | 97.9 |

specific hazards model for $\lambda(t,v|z)$ and a proportional hazards model for $\lambda(t|z) = \lambda_0(t)\exp(\beta z)$, with the mark-specific vaccine efficacy $\mathrm{VE}(v) = 1 - 2v$ and the marginal $\mathrm{VE} = 1 - \exp(\beta) = 0$. The standard Wald test, denoted by $T_w$, under the marginal Cox model is often used to test for the vaccine efficacy. As expected, the standard Wald test shows no power (Table 3). It is incapable of revealing any vaccine efficacy or that the vaccine efficacy depends on the mark, thus missing the important scientific finding that the vaccine protects against viruses with smaller mark values ($V < 0.5$) and increases risk of infection with viruses with larger mark values ($V > 0.5$). The example we constructed here shows the weakness of using the standard approach that ignores the mark and is what motivates the present research.

**4. Application.** The first preventive HIV vaccine efficacy trial was carried out in North America and The Netherlands, and enrolled 5403 HIV-

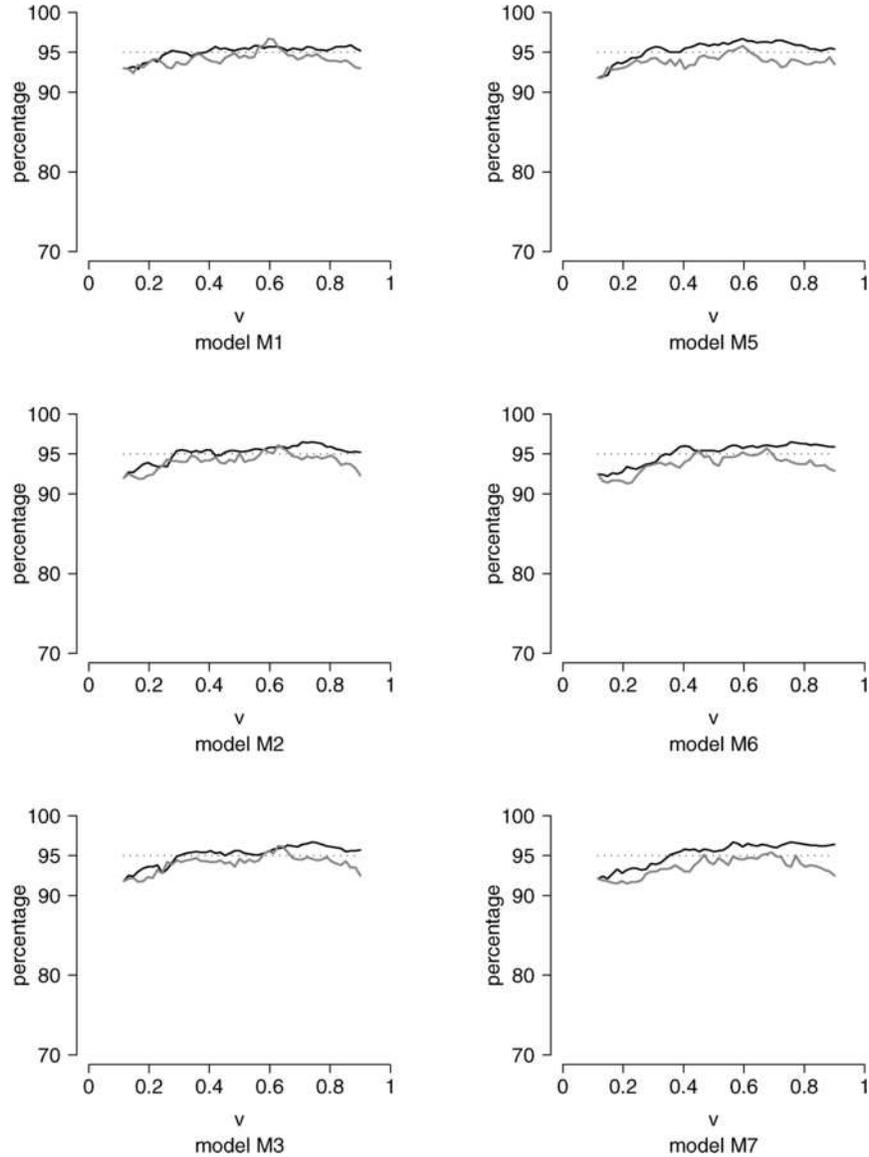

Fig. 3. *Plots of the pointwise coverage probabilities for* VE($v$) *(gray lines) and for* CV($v$) *(solid lines), based on* $n = 500$, $h = 0.1$ *and 1000 repetitions. The models on the left panel are M1, M2 and M3. The models on the right panel are M5, M6 and M7.*

negative volunteers at risk for acquiring HIV infection [4]. Volunteers were randomized in a 2:1 ratio to receive a recombinant glycoprotein 120 vaccine (AIDSVAX) or placebo, and were monitored for HIV infection at semi-annual HIV testing visits for 36 months. The primary objective was to assess

TABLE 3
*Comparison of the standard Wald test with the proposed tests $T_a^{(1)}$, $T_{m1}^{(1)}$, $T_{m2}^{(1)}$, $T_a^{(2)}$, $T_{m1}^{(2)}$ and $T_{m2}^{(2)}$ at the nominal level 0.05*

| | | Power | | | | | | |
|---|---|---|---|---|---|---|---|---|
| $n$ | $h$ | $T_w$ | $T_a^{(1)}$ | $T_{m1}^{(1)}$ | $T_{m2}^{(1)}$ | $T_a^{(2)}$ | $T_{m1}^{(2)}$ | $T_{m2}^{(2)}$ |
| 500 | 0.05 | 5.9 | 14.9 | 24.2 | 16.6 | 98.0 | 99.4 | 97.3 |
| | 0.1 | — | 23.9 | 35.7 | 16.0 | 99.6 | 100 | 99.8 |
| | 0.15 | — | 27.9 | 39.1 | 15.7 | 99.9 | 100 | 99.9 |
| 800 | 0.05 | 6.1 | 32.4 | 39.6 | 15.0 | 100 | 100 | 99.6 |
| | 0.1 | — | 43.1 | 51.5 | 13.8 | 100 | 100 | 100 |
| | 0.15 | — | 46.0 | 53.3 | 13.9 | 100 | 100 | 100 |

VE using the standard Cox model, and a secondary objective was to test $H_{10}: \mathrm{VE}(t,v) = 0$ and $H_{20}: \mathrm{VE}(t,v) = \mathrm{VE}(t)$ for three different mark variables $V$ defined in terms of the percent mismatch of aligned amino acid sequences (for each infecting HIV sequence compared to the HIV sequence [named GNE8] contained in the AIDSVAX construct) in three subregions of HIV-gp120. For brevity, in this article we consider only one mark $V$, defined as the percent mismatch of amino acids in the whole gp120 region (581 amino acids long), where all possible mismatches of particular pairs of amino acids (e.g., A versus C) are weighted by the estimated probability of interchange [12]. The distance is based on the gp120 region because this region contains neutralizing epitopes that potentially can induce anti-HIV antibody responses that prevent HIV infection [21]; the vaccine was designed to protect by stimulating high titer antibodies that neutralize exposing HIVs. Of the 368 individuals infected during the trial, 32 had missing marks. Of the remaining 336 samples, all marks were unique (217 vaccine; 119 placebo).

The vaccine efficacy is estimated and tested by adjusting for two covariates: age (ranging 18–62 years with mean of 36.5) and behavioral risk score (taking values 0–7) as defined in [4]. It is relevant to adjust for these covariates because they predict infection rate and because trial participants with different values of these covariates may be exposed to HIV strains with different distributions of $V$. Both covariates are considered as continuous variables. The histograms of the rescaled mark values, ages in years and behavioral risk scores are plotted in Figure 4. We denote the treatment indicator by $z_1$ ($z_1 = 1$ for the vaccine and $z_1 = 0$ for the placebo), age by $z_2$ and behavioral risk score by $z_3$, and denote the corresponding coefficient functions by $\beta_1(v)$, $\beta_2(v)$ and $\beta_3(v)$. Fitting model (2) with $h = 0.3$, the plots of the estimates for $\beta_1(v)$, $\beta_2(v)$ and $\beta_3(v)$ and their pointwise confidence bands are given in Figure 5. The plots of $\widehat{\mathrm{VE}}(v)$ and $\widehat{\mathrm{CV}}(v)$ with their corresponding pointwise confidence bands adjusting for the two covariates $z_2$ and $z_3$ are given in Figure 6.

Adjusting for age and behavioral risk score, the Wald test statistic for testing the marginal $\text{VE} = 0$ using the standard Cox model is $-0.978$, yielding a $p$-value of 0.328 for the two-sided alternative and 0.164 for the monotone alternative. Our test with the test statistic $T_a^{(1)}$ for $H_{10} : \text{VE}(v) = 0$ for all $v$ versus the general alternative $H_{1a}$ yields a $p$-value of 0.1532. The $p$-values for testing against the monotone alternative $H_{1m}$ are 0.0916 for $T_{m1}^{(1)}$ and 0.0228 for $T_{m2}^{(1)}$. These results give some, albeit weak, evidence of nonzero vaccine efficacy for at least one mark value; see Figure 6.

In addition, adjusting for age and behavioral risk score, we conducted the tests to evaluate whether the vaccine efficacy varies with the mark. The $p$-value for testing $H_{20}$ that $\text{VE}(v)$ does not depend on $v$ versus the general alternative $H_{2a}$ is 0.2067 for the test statistic $T_a^{(2)}$. The $p$-value for testing for the monotone alternative $H_{2m}$ is 0.9363 for the test statistic $T_{m1}^{(2)}$ and 0.9047 for the test statistic $T_{m2}^{(2)}$. These $p$-values are expected given the plots in Figure 6 where $\widehat{\text{VE}}(v)$ shows some tendency to increase with $v$.

**5. Discussion.** This article developed inference techniques for the proportional hazards model with a continuous mark variable, including nonparametric methods for estimation and testing of mark-specific regression functions. These techniques can be used to estimate mark-specific vaccine efficacy ($\text{VE}(v)$) and cumulative mark-specific vaccine efficacy ($\text{CV}(v)$) with simultaneous confidence bands, and to test hypotheses for $\text{VE}(v)$, while adjusting for time-dependent covariate effects. The testing procedures based on the statistics $T_{m1}^{(1)}$ and $T_{m2}^{(2)}$ showed greatest power in simulations and are recommended for testing $\text{VE}(v) = 0$ for all $v$ and for testing $\text{VE}(v)$ independent of $v$, respectively.

An alternative approach to the continuous mark-specific PH model would be a similar model that treats the mark variable as ordinal categorical. We focused on a continuous mark because (i) it most naturally suits the HIV vaccine application, as the choice of $K$ bins for categorizing the marks would be arbitrary and (ii) testing $\beta(v) = \beta$ can often be done with greater power than testing equality of the cause-specific regression coefficients $\beta_1 = \cdots = \beta_K$.

As is well known for a discrete mark-specific hazard function, the interpretation of the continuous mark-specific hazard function $\lambda(t, v)$ is restricted to actual study conditions, that is, it is the instantaneous rate of failure in the presence of all of the circulating competing risks (i.e., is a "crude" hazard in the terminology of Prentice et al. [14]). However, often the main scientific interest is in the "net" mark-specific hazard, the instantaneous rate of failure by mark $v$ in the absence of any other competing risks, but unfortunately this parameter is not identified except under untestable assumptions such as mutual independence of all of the notional (latent) mark-specific failure times [18]. This problem necessitates careful interpretation of inferences in the mark-specific PH model.

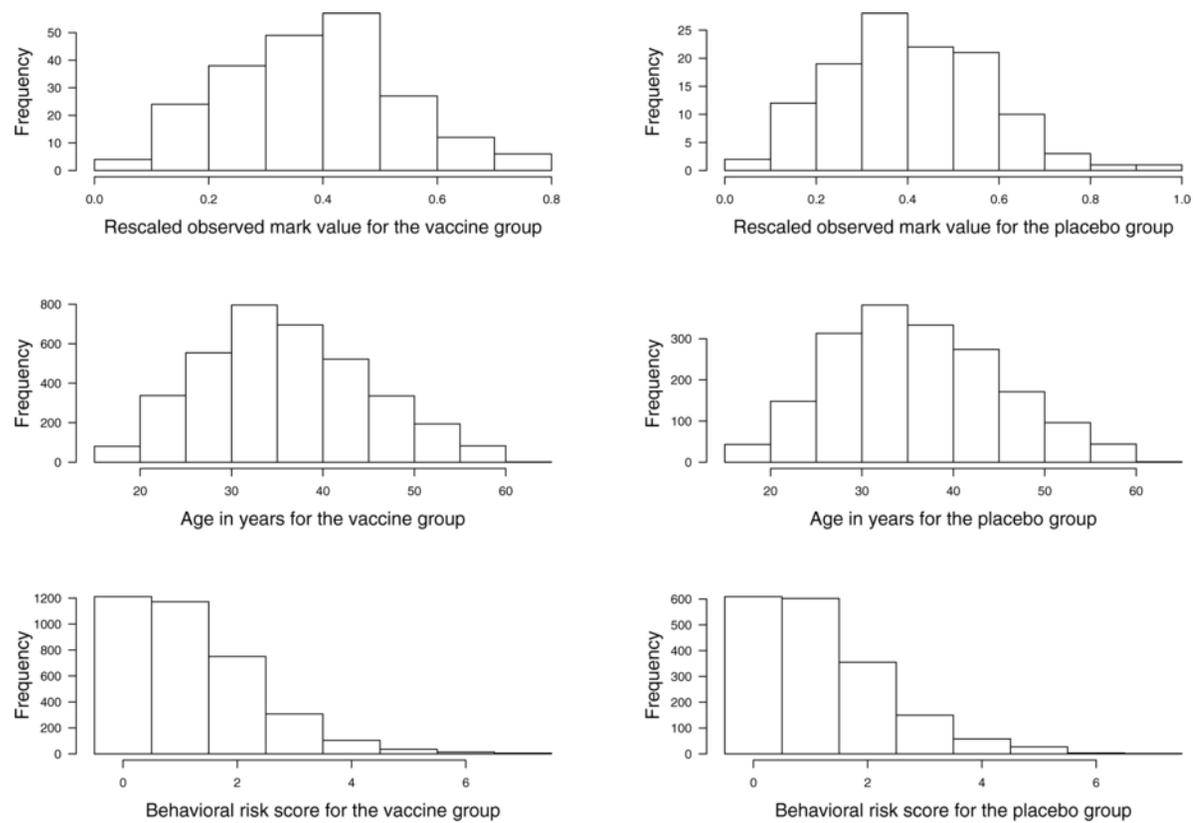

Fig. 4. *Histograms for the observed mark values, ages in years and behavioral risk scores. The left panel is for the vaccine group and the right panel is for the placebo group.*

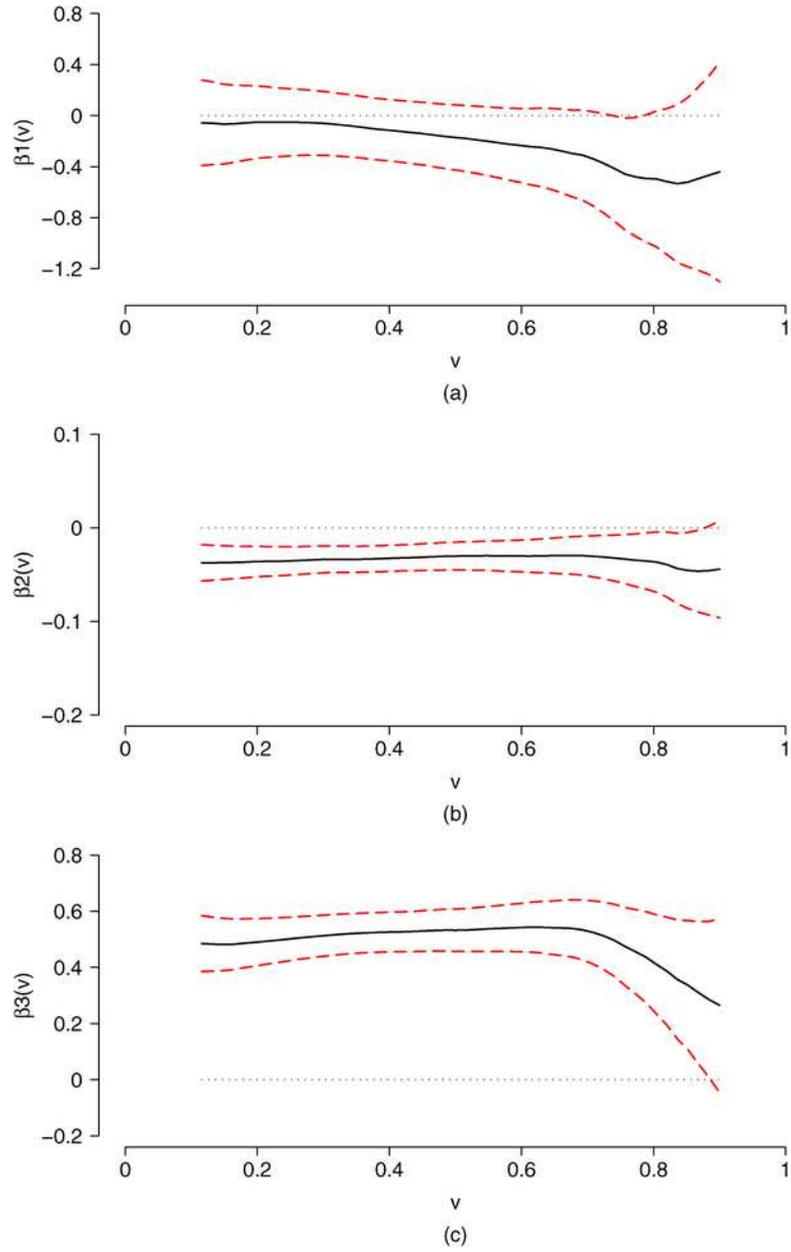

FIG. 5. *Plots of the estimated regression coefficients $\beta_1(v)$, $\beta_2(v)$ and $\beta_3(v)$ and their 95% pointwise confidence bands for the vaccine trial data with $h = 0.3$.*

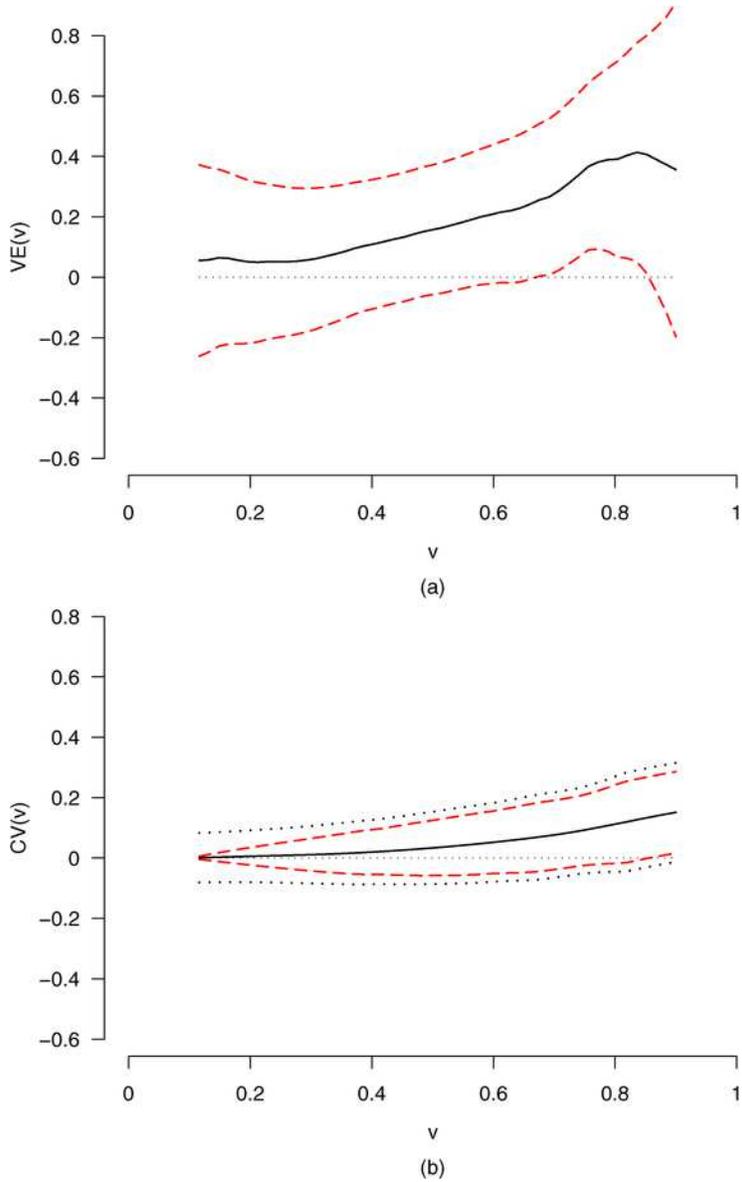

FIG. 6. *Plots of the estimates of* VE(v) *and* CV(v) *and their confidence bands for the vaccine trial data with* $h = 0.3$. *The dashed lines are 95% pointwise confidence bands and the dotted lines are 95% simultaneous confidence bands.*

For the HIV vaccine trial example, the crude mark-specific hazard can be factored as

$$\lambda(t,v|z) = \lambda_E(t,v|z) \times \lambda_{PC}(t|v,z) \tag{14}$$

where $\lambda_E(t,v|z)$ is the intensity of exposure to strain $v$ for participants with covariates $z$ and $\lambda_{PC}(t|v,z)$ (the "per-contact" transmission hazard) is the same as $\lambda(t,v|z)$ except that it further conditions on the (unobserved) presence of exposure to a virus with genetic distance $v$ during $[t, t+dt)$. Exposure can arise from unprotected sex or sharing a needle with an individual infected with strain $v$. Therefore the identified parameter measures a mixture of vaccine/placebo-group differences in mark-specific exposure rates and in conditional mark-specific per-exposure transmission probabilities, whereas biological interest is in

$$\text{VE}^{PC}(t|v,z_2) = 1 - \frac{\lambda_{PC}(t|v,1,z_2)}{\lambda_{PC}(t|v,0,z_2)}$$

as a measure of vaccine efficacy. However, as data are not available for estimating the relative intensity $\lambda_E(t,v|1,z_2)/\lambda_E(t,v|0,z_2)$, our approach is to use

$$\text{VE}(t,v|z_2) = 1 - \frac{\lambda(t,v|1,z_2)}{\lambda(t,v|0,z_2)}$$

as the target estimand, and assume identical exposure rates between the two groups, so this target has the same interpretation as $\text{VE}^{PC}(t|v,z_2)$. Reliance on this assumption demonstrates the value of including covariates $z_2$ that predict mark-specific exposure into the mark-specific PH model: the richer the covariate information the more likely $\text{VE}(t,v|z_2)$ reflects biological vaccine efficacy. Gilbert, McKeague and Sun [5] provided further discussion of the interpretation of mark-specific hazard ratios.

The usefulness of our approach relies on the validity of the mark-specific proportional hazards model. Lin, Wei and Ying [9] developed goodness-of-fit tests for the standard Cox model based on martingale residuals, and their tests can be extended to the present setting by using the mark-specific martingale residuals

$$\hat{M}_i(t,v) = \int_0^t \int_a^v [N_i(ds,du) - Y_i(s)\exp((\hat{\beta}(u))^T Z_i)\hat{\Lambda}_0(ds,du)], \tag{15}$$

for $i = 1, \ldots, n$. These residuals may be interpreted as the difference at time $t$ between the observed and the predicted number of events with mark less than $v$ for the $i$th subject, and are informative about model misspecification. It can be checked that $n^{-1/2}\sum_{i=1}^n \hat{M}_i(t,v) = o_p(1)$. This property is similar to that in the standard Cox model, where the sum of the martingale residuals is exactly zero. The difference here is caused by the kernel smoothing in a neighborhood of $v$. Because $\beta(v)$ is treated nonparametrically, the checking of the model (2) needs further development and has additional issues related to the bandwidth. This would need a thorough treatment that is beyond the scope of the present paper.

Finally, we caution that the method proposed here requires large sample sizes to work well as demonstrated in the simulation study. This is the result of $\beta(v)$ being treated nonparametrically: the estimation of $\beta(v)$ utilizes only the observed failures with marks in a neighborhood of $v$. Although this does not cause a problem in our application to the first preventive HIV vaccine trial (which has a sample size of 5403), one needs to be careful in applying the method to situations with small sample sizes.

## APPENDIX

The following lemma is an extension of Theorem 5.7 of van der Vaart [19] and will be used to prove the uniform consistency of $\hat{\beta}(v)$.

LEMMA A.1. *Let $Q_n(v, \theta)$ be random functions and let $Q(v, \theta)$ be a fixed function of $(v, \theta) \in [a, b] \times \Theta$, $\Theta \subset \mathbb{R}^p$. Let $\beta(v)$ be a fixed function of $v \in [a, b]$ taking values in $\Theta$. Assume that $\sup_{v,\theta} |Q_n(v, \theta) - Q(v, \theta)| \xrightarrow{P} 0$ and that for every $\varepsilon > 0$ there exists an $\eta > 0$ such that $\sup_{\|\theta - \beta(v)\| > \varepsilon} Q(v, \theta) < Q(v, \beta(v)) - \eta$ for $v \in [a, b]$. Then for any sequence of estimators $\hat{\beta}(v)$, with $Q_n(v, \hat{\beta}(v)) > Q_n(v, \beta(v)) - o_p(1)$ uniformly in $v \in [a, b]$, we have $\hat{\beta}(v) \xrightarrow{P} \beta(v)$ uniformly in $v \in [a, b]$.*

PROOF. For every $\varepsilon > 0$, there exists an $\eta > 0$ such that

$$\left\{\sup_v \|\hat{\beta}(v) - \beta(v)\| > \varepsilon\right\} \subset \bigcup_v \{\|\hat{\beta}(v) - \beta(v)\| > \varepsilon\}$$
$$\subset \bigcup_v \{Q(v, \hat{\beta}(v)) < Q(v, \beta(v)) - \eta\}.$$

Since $Q_n(v, \hat{\beta}(v)) > Q_n(v, \beta(v)) - o_p(1) \xrightarrow{P} Q(v, \beta(v))$, uniformly in $v \in [a, b]$, we have $Q_n(v, \hat{\beta}(v)) > Q(v, \beta(v)) - o_p(1)$, uniformly in $v \in [a, b]$. It follows that

$$\bigcup_v \{Q(v, \hat{\beta}(v)) < Q(v, \beta(v)) - \eta\}$$
$$\subset \bigcup_v \{Q(v, \hat{\beta}(v)) < Q_n(v, \hat{\beta}(v)) - \eta + o_p(1)\}$$
$$= \left\{\inf_v (Q(v, \hat{\beta}(v)) - Q_n(v, \hat{\beta}(v))) < -\eta + o_p(1)\right\}$$

$$= \left\{ \sup_v (Q_n(v, \hat{\beta}(v)) - Q(v, \hat{\beta}(v))) > \eta - o_p(1) \right\}$$

$$\subset \left\{ \sup_v |Q_n(v, \hat{\beta}(v)) - Q(v, \hat{\beta}(v))| > \eta - o_p(1) \right\},$$

whose probability goes to 0 by the uniform convergence of $Q_n(v, \theta)$ to $Q(v, \theta)$. Hence $P\{\sup_v \|\hat{\beta}(v) - \beta(v)\| > \varepsilon\} \to 0$. □

The following lemma is used to prove Theorems 3 and 4. Let $N = \sum_{i=1}^n N_i$ and $M = \sum_{i=1}^n M_i$.

LEMMA A.2. *Under conditions* (A.1)–(A.4), $n^{-1}N(t,v) \xrightarrow{P} EN_i(t,v)$, *uniformly in* $(t,u) \in [0,\tau] \times [0,1]$, *and* $n^{-1/2}M(t,v)$ *converges weakly to a mean-zero continuous Gaussian random field* $G(t,v)$, $(t,v) \in [0,\tau] \times [0,1]$, *with independent increments and* $\mathrm{Var}(G(t,v)) = \int_0^t \int_0^v \lambda_0(s,u) s^{(0)}(s, \beta(u))\, ds\, du$.

PROOF. We treat $\omega_i = (X_i, \delta_i, V_i)$, $i = 1, \ldots, n$, as a random sample from a probability distribution $P$ on a measurable space $(\mathcal{X}, \mathcal{A})$, with $\mathcal{X} = [0, \infty) \times \{0,1\} \times [0,1]$ and $\mathcal{A}$ its Borel $\sigma$-field. Let $\mathcal{F}$ be the class of all indicator functions $f_{t,v}: \mathcal{X} \longrightarrow R$, where $f_{t,v}(\omega_i) = I([0,t] \times \{1\} \times [0,v])(\omega_i) = I(X_i \leq t, \delta_i = 1, V_i \leq v)$, for $0 \leq t \leq \tau, 0 \leq v \leq 1$. Then $n^{-1}N(t,v) = n^{-1}\sum_{i=1}^n f_{t,v}(\omega_i)$. Let $\|f_{t,v}\|_{P,r} = (P|f_{t,v}|^r)^{1/r}$ be $L_r(P)$-norm of $f_{t,v}$.

Let $0 = t_0 < t_1 < \cdots < t_K = \tau$ and $0 = v_0 < v_1 < \cdots < v_J = 1$ be partitions of the intervals $[0,\tau]$ and $[0,1]$. Define the bracketing functions $l_{kj} = N_i(t_{k-1}, v_{j-1})$ and $u_{kj} = N_i(t_k, v_j)$, for $k = 1, \ldots, K$, $j = 1, \ldots, J$. Then for any $f_{t,v} \in \mathcal{F}$, there is a bracket $[l_{kj}, u_{kj}]$ such that $f_{t,v} \in [l_{kj}, u_{kj}]$. And

$$\|u_{kj} - l_{kj}\|_{P,1} \leq E(N_i(t_k, v_j) - N_i(t_{k-1}, v_{j-1}))$$

$$= \int_0^{t_k} \int_0^{v_j} \lambda_0(s,x) s^{(0)}(s, \beta(x))\, ds\, dx$$

$$- \int_0^{t_{k-1}} \int_0^{v_{j-1}} \lambda_0(s,x) s^{(0)}(s, \beta(x))\, ds\, dx$$

$$\leq \int_{t_{k-1}}^{t_k} \int_0^1 \lambda_0(s,x) s^{(0)}(s, \beta(x))\, ds\, dx$$

$$+ \int_0^\tau \int_{v_{j-1}}^{v_j} \lambda_0(s,x) s^{(0)}(s, \beta(x))\, ds\, dx$$

$$\leq C_1(t_k - t_{k-1}) + C_2(v_j - v_{j-1}),$$

where $C_1$ and $C_2$ are some positive constants. For any $\varepsilon > 0$, choose the grid points such that $t_k - t_{k-1} < \varepsilon$ and $v_j - v_{j-1} < \varepsilon$. Then $\|u_{kj} - l_{kj}\|_{P,1} \leq [C_1 + C_2]\varepsilon$. Hence, the bracketing number $N_{[\cdot]}(\varepsilon, \mathcal{F}, L_1(P))$ is of the polynomial

order $(1/\varepsilon)^2$. By the Glivenko–Cantelli theorem (Theorem 19.4 of van der Vaart [19]), $n^{-1}N(t,v) \xrightarrow{P} EN_i(t,v)$, uniformly in $(t,v) \in [0,\tau] \times [0,1]$.

Next, consider the processes $\{M_i(t,v), 0 \leq t \leq \tau, 0 \leq v \leq 1\}$, $i = 1, \ldots, n$, as a random sample from a probability distribution $P$ on a measurable space $(\mathcal{X}, \mathcal{A})$. Let $\mathcal{F}$ be the class of coordinate projections $f_{t,v} : \mathcal{X} \longrightarrow R$, where $f_{t,v}(M_i) = M_i(t,v)$, for $0 \leq t \leq \tau, 0 \leq v \leq 1$. The process $\{M_i(t,v), 0 \leq t \leq \tau, 0 \leq v \leq 1\}$ is determined by the $\{X_i, \delta_i, \delta_i V_i, Z_i\}$.

Again, let $0 = t_0 < t_1 < \cdots < t_K = \tau$ and $0 = v_0 < v_1 < \cdots < v_J = 1$ be the partitions of the intervals $[0, \tau]$ and $[0, 1]$. Define the bracketing functions $l_{kj} = N_i(t_{k-1}, v_{j-1}) - \int_0^{t_k} \int_0^{v_j} Y_i(s)\lambda(s, x|Z_i(s))\,ds\,dx$ and $u_{kj} = N_i(t_k, v_j) - \int_0^{t_{k-1}} \int_0^{v_{j-1}} Y_i(s)\lambda(s, x|Z_i(s))\,ds\,dx$, for $k = 1, \ldots, K$, $j = 1, \ldots, J$. Then for any $f_{t,v} \in \mathcal{F}$, there is a bracket $[l_{kj}, u_{kj}]$ such that $f_{t,v} \in [l_{kj}, u_{kj}]$. The bracket size is

$$\|u_{kj} - l_{kj}\|_{P,2} \leq \|N_i(t_k, v_j) - N_i(t_{k-1}, v_{j-1})\|_{P,2}$$
$$+ \left\|\int_0^{t_k} \int_0^{v_j} Y_i(s)\lambda(s, x|Z_i(s))\,ds\,dx \right.$$
$$\left. - \int_0^{t_{k-1}} \int_0^{v_{j-1}} Y_i(s)\lambda(s, x|Z_i(s))\,ds\,dx\right\|_{P,2}$$
$$\leq [C_1(t_k - t_{k-1}) + C_2(v_j - v_{j-1})]^{1/2},$$

where $C_1$ and $C_2$ are some positive constants. For any $\varepsilon > 0$, choose the grid points such that $t_k - t_{k-1} < \varepsilon$ and $v_j - v_{j-1} < \varepsilon$. Then $\|u_{kj} - l_{kj}\|_{P,2} \leq [C_1 + C_2]^{1/2}\varepsilon^{1/2}$. Hence, the bracketing number $N_{[\cdot]}(\varepsilon^{1/2}, \mathcal{F}, L_2(P))$ is of the polynomial order $(1/\varepsilon)^2$. Thus, $N_{[\cdot]}(\varepsilon, \mathcal{F}, L_2(P))$ is of the polynomial order $(1/\varepsilon)^4$. So the bracketing integral $J_{[\cdot]}(1, \mathcal{F}, L_2(P)) < \infty$. By the Donsker theorem (Theorem 19.5 of van der Vaart [19]), $n^{-1/2}M = \{n^{-1/2}\sum_{i=1}^n M_i(t,v), 0 \leq t \leq \tau, 0 \leq v \leq 1\}$ converges weakly to a mean-zero Gaussian process $G(t,v)$, $(t,v) \in [0,\tau] \times [0,1]$, which can be constructed to have continuous paths by Theorem 18.14 and Lemma 18.15 of van der Vaart [19].

Now we show that $G(t,v)$ has independent increments. Note that for $t_1 \leq t_2$ and $v_1 \leq v_2$, the covariance of $G(t_1, v_1)$ and $G(t_2, v_2) - G(t_1, v_1)$ is $E\{M_i(t_1, v_1) \times (M_i(t_2, v_2) - M_i(t_1, v_1))\}$. By Aalan and Johansen [1], $M_i(t, v_1)$ and $M_i(t, v_2) - M_i(t, v_1)$, $0 \leq t \leq \tau$, are orthogonal square integrable martingales for $0 \leq v_1 \leq v_2 \leq 1$. It follows that

$$E\{M_i(t_1, v_1)(M_i(t_2, v_2) - M_i(t_1, v_1))\}$$
$$= E\{M_i(t_1, v_1)(M_i(t_2, v_2) - M_i(t_2, v_1))\}$$
$$\quad + E\{M_i(t_1, v_1)(M_i(t_2, v_1) - M_i(t_1, v_1))\}$$
$$= 0.$$

Hence $G(t_1, v_1)$ and $G(t_2, v_2) - G(t_1, v_1)$ are independent. □

PROOF OF THEOREM 1. It is easy to check that the conditions of Lemma 1 of Sun and Wu [17] are satisfied under Condition A. It follows that $\tilde{W}_A(v)$ converges weakly to a vector of continuous mean-zero Gaussian random processes, $W_A(v)$, $v \in [a, b]$. Now we show that $W_A(v)$ has independent increments. Let $w_i(t, v) = \int_a^v \int_0^t A(u)[Z_i(t) - s^{(1)}(t, \beta(u))/s^{(0)}(t, \beta(u))]M_i(dt, du)$. Then $\tilde{W}_A(v) = n^{-1/2} \sum_{i=1}^n w_i(\tau, v)$. For $a \leq v_1 \leq v_2 \leq b$, the covariance matrix of $W_A(v_1)$ and $W_A(v_2) - W_A(v_1)$ is equal to $E\{w_i(\tau, v_1)(w_i(\tau, v_2) - w_i(\tau, v_1))^T\}$. Since $M_i(t, v_1)$ and $M_i(t, v_2) - M_i(t, v_1)$, $0 \leq t \leq \tau$, are orthogonal square integrable martingales, it follows that $w_i(t, v_1)$ and $w_i(t, v_2) - w_i(t, v_1)$, $0 \leq t \leq \tau$, are orthogonal square integrable martingales. Hence $E\{w_i(\tau, v_1)(w_i(\tau, v_2) - w_i(\tau, v_1))^T\} = 0$. So $W_A(v)$, $v \in [a, b]$, is a vector of mean-zero Gaussian random processes with independent increments.

Further, the covariance matrix of $W_A(v)$ is

$$E\{w_i(\tau, v)(w_i(\tau, v))^T\}$$
$$= E\left\{\int_a^v \int_0^\tau A(u)\left[Z_i(t) - \frac{s^{(1)}(t, \beta(u))}{s^{(0)}(t, \beta(u))}\right]^{\otimes 2} A(u)\, N_i(dt, du)\right\}$$
$$= E\left\{\int_a^v \int_0^\tau A(u)\left[Z_i(t) - \frac{s^{(1)}(t, \beta(u))}{s^{(0)}(t, \beta(u))}\right]^{\otimes 2}\right.$$
$$\left.\times A(u) y(t|Z_i(t))\lambda(t, u|Z_i(t))\, dt\, du\right\}$$
$$= \int_a^v A(u) E\left\{\int_0^\tau \left[Z_i(t) - \frac{s^{(1)}(t, \beta(u))}{s^{(0)}(t, \beta(u))}\right]^{\otimes 2} y(t|Z_i(t))\lambda(t, u|Z_i(t))\, dt\right\}$$
$$\times A(u)\, du$$
$$= \int_a^v A(u)\Sigma(u)A(u)\, du.$$

This completes the proof of Theorem 1. □

PROOF OF THEOREM 2. We shall prove Theorem 2 by verifying the conditions of Lemma A.1.

Let

$$\eta_n(u, \theta) = n^{-1} \sum_{i=1}^n \int_0^u \int_0^\tau [\theta^T Z_i(t) - \log(S^{(0)}(t, \theta))]\, N_i(dt, du),$$

$$\xi_n(u, \theta) = n^{-1} \sum_{i=1}^n \int_0^u \int_0^\tau [\theta^T Z_i(t) - \log(s^{(0)}(t, \theta))]\, N_i(dt, du),$$

$$Q_n(v, \theta) = n^{-1} l(v, \theta) + n^{-1} \log n \int_0^1 K_h(u - v) N(\tau, du).$$

Then by Condition A, $\eta_n(v,\theta) = \xi_n(v,\theta) + O_p(n^{-1/2})$ and

$$Q_n(v,\theta) = \int_0^1 K_h(u-v)\eta_n(du,\theta)$$
$$= \int_0^1 K_h(u-v)\xi_n(du,\theta) + O_p(n^{-1/2}h^{-1}),$$

uniformly in $(v,\theta) \in [0,1] \times [-M,M]$, for $M > 0$. By application of the Glivenko–Cantelli and Donsker theorems, similarly to the proofs of Lemma A.2 and Theorem 1, $\xi_n(v,\theta) = \xi(v,\theta) + O_p(n^{-1/2})$, uniformly in $(v,\theta) \in [0,1] \times [-M,M]$, with

$$\xi(v,\theta) = E\left[\int_0^u \int_0^\tau [\theta^T Z_i(t) - \log(s^{(0)}(t,\theta))]N_i(dt,du)\right].$$

It follows that $Q_n(v,\theta) = Q(v,\theta) + O_p(n^{-1/2}h^{-1})$, uniformly in $(v,\theta) \in [a,b] \times [-M,M]$, where

$$Q(v,\theta) = E\left[\int_0^\tau [\theta^T Z_i(t) - \log(s^{(0)}(t,\theta))]\lambda_0(t,v)\exp(\beta^T(v)Z_i(t))Y_i(t)\,dt\right].$$

Now we show that $\beta(v)$ is the well-separated point of maximum of $Q(v,\theta)$ for $v \in [0,1]$. Note that

$$\partial Q(v,\theta)/\partial \theta = E\left[\int_0^\tau \left[Z_i(t) - \frac{s^{(1)}(t,\theta)}{s^{(0)}(t,\theta)}\right]\lambda_0(t,v)\exp(\beta^T(v)Z_i(t))Y_i(t)\,dt\right]$$

$$\partial^2 Q(v,\theta)/\partial \theta^2 = -E\left[\int_0^\tau \left\{\frac{s^{(2)}(t,\theta)}{s^{(0)}(t,\theta)} - \left(\frac{s^{(1)}(t,\theta)}{s^{(0)}(t,\theta)}\right)^{\otimes 2}\right\}\right.$$
$$\left. \times \lambda_0(t,v)\exp(\beta^T(v)Z_i(t))Y_i(t)\,dt\right].$$

We have $\partial Q(v,\beta(v))/\partial \theta = 0$, and for every $\varepsilon > 0$ there exists an $\eta > 0$ such that $\sup_{\|\theta-\beta(v)\|>\varepsilon} Q(v,\theta) < Q(v,\beta(v)) - \eta$ for $v \in [a,b]$, under condition (A.3), by Taylor expansion and continuity. Further, since $Q_n(v,\theta) \xrightarrow{P} Q(v,\theta)$, $\partial Q_n(v,\theta)/\partial \theta \xrightarrow{P} \partial Q(v,\theta)/\partial \theta$, and $\partial^2 Q_n(v,\theta)/\partial \theta^2 \xrightarrow{P} \partial^2 Q(v,\theta)/\partial \theta^2$ uniformly in $(v,\theta) \in [a,b] \times [-M,M]$, and $-\tilde{M} < \beta(v) < \tilde{M}$ for $a \leq v \leq b$ for some $\tilde{M} < M$, we have for every $\alpha > 0$ there exists an $n_0$ such that $P(-M \leq \hat{\beta}(v) \leq M, a \leq v \leq b) > 1 - \alpha$ for $n \geq n_0$.

Therefore, for every $\varepsilon > 0$,

$$P\left(\sup_{a \leq v \leq b} \|\hat{\beta}(v) - \beta(v)\| > \varepsilon\right)$$
$$\leq \alpha + P\left(\sup_{a \leq v \leq b} \|\hat{\beta}(v) - \beta(v)\| > \varepsilon, -M \leq \hat{\beta}(v) \leq M, a \leq v \leq b\right)$$
$$\to \alpha$$

as $n \to \infty$, by the previous checking of the conditions of Lemma A.1 together with $Q_n(v, \hat{\beta}(v)) \geq Q_n(v, \beta(v))$. Since $\alpha$ is arbitrary, we have $P(\sup_{a \leq v \leq b} \|\hat{\beta}(v) - \beta(v)\| > \varepsilon) \to 0$. $\square$

PROOF OF THEOREM 3. In the proof of this theorem, we set $\beta = \beta(v)$ for simplicity. Note that under Condition A, using a second-order Taylor expansion for $\lambda(t, u|Z_i(t))$ in the neighborhood of $v$, we have

$$n^{-1/2} \left| \sum_{i=1}^{n} \int_0^1 \int_0^{\tau} K_h(u-v) \left[ Z_i(t) - \frac{S^{(1)}(t,\beta)}{S^{(0)}(t,\beta)} \right] Y_i(t) \right.$$

$$\left. \times [\lambda(t, v|Z_i(t)) - \lambda(t, u|Z_i(t))] \, dt \, du \right|$$

$$= O_p(n^{1/2} h^2),$$

uniformly in $v \in [0, 1]$. It follows that

$$n^{-1/2} U(v, \beta) = n^{-1/2} \sum_{i=1}^{n} \int_0^1 \int_0^{\tau} K_h(u-v) \left[ Z_i(t) - \frac{S^{(1)}(t,\beta)}{S^{(0)}(t,\beta)} \right]$$

$$\times [N_i(dt, du) - Y_i(t) \lambda(t, v|Z_i(t)) \, dt \, du]$$

$$= n^{-1/2} \sum_{i=1}^{n} \int_0^1 \int_0^{\tau} K_h(u-v) \left[ Z_i(t) - \frac{S^{(1)}(t,\beta)}{S^{(0)}(t,\beta)} \right] M_i(dt, du)$$

$$+ O_p(n^{1/2} h^2),$$

uniformly in $v \in [0, 1]$.

Next, we show that for each $v$, $n^{-1/2} h^{1/2} U(v, \beta)$ converges weakly to a normal distribution. By Lemma A.2, $n^{-1/2} M(t, v)$ converges weakly to a mean-zero Gaussian process. By Condition A, $\|S^{(j)}(t, \beta) - s^{(j)}(t, \beta)\| = o_p(n^{-1/2+\delta})$, uniformly in $t$ for $j = 0, 1$, for $0 < \delta < 1/2$. Note that $n^{-1/2+\delta} h^{-1/2} = o(1)$ for $\delta = 1/4$ as $nh^2 \to \infty$. We have $h^{1/2} K_h(u-v) \|S^{(j)}(t, \beta) - s^{(j)}(t, \beta)\|$ goes in probability to zero. Applying Lemma 2 of Gilbert, McKeague and Sun [5], we have

$$n^{-1/2} h^{1/2} U(\beta(v))$$

$$= n^{-1/2} h^{1/2} \sum_{i=1}^{n} \int_0^1 \int_0^{\tau} K_h(u-v) \left[ Z_i(t) - \frac{s^{(1)}(t,\beta)}{s^{(0)}(t,\beta)} \right]$$

$$\times M_i(dt, du) + O_p(n^{1/2} h^{5/2}) + o_p(1)$$

(16)

$$= n^{-1/2} h^{1/2} \sum_{i=1}^{n} \int_0^1 \int_0^{\tau} K_h(u-v) \left[ Z_i(t) - \frac{s^{(1)}(t, \beta(u))}{s^{(0)}(t, \beta(u))} \right]$$

$$\times M_i(dt, du) + O_p(n^{1/2}h^{5/2}) + o_p(1)$$

$$= h^{1/2} \int_0^1 K_h(u - v)\tilde{W}_I(du) + O_p(n^{1/2}h^{5/2}) + o_p(1),$$

where $\tilde{W}_I(v)$ is defined in (6) with $A = I$ and $a = 0$.

Since $\tilde{W}_I(v) \xrightarrow{\mathcal{D}} W_I(v)$ by Theorem 1, by the almost sure representation theorem ([16], page 47), there exist $\tilde{W}_I^*(v)$ and $W_I^*(v)$ on some probability space that have the same distributions and sample paths as $\tilde{W}_I(v)$ and $W_I(v)$, respectively, such that $\tilde{W}_I^*(v) \xrightarrow{\text{a.s.}} W_I^*(v)$ uniformly in $v \in [0, 1]$. Hence $\int_0^1 K_h(u-v)\tilde{W}_I^*(du) = \int_0^1 K_h(u-v)W_I^*(du) + O_p(n^{-1/2}h^{-1})$ by integration by parts since $K(\cdot)$ has bounded variation. It follows that

$$h^{1/2}\int_0^1 K_h(u-v)\tilde{W}_I(du) \stackrel{\mathcal{D}}{=} h^{1/2}\int_0^1 K_h(u-v)\tilde{W}_I^*(du)$$

$$= h^{1/2}\int_0^1 K_h(u-v)W_I^*(du) + O_p(n^{-1/2}h^{-1/2}).$$

Since $W_I^*(v)$ is a Gaussian martingale with covariance matrix of $\int_0^v \Sigma(u)\,du$, and $h^{1/2}\int_0^1 K_h(u-v)W_I^*(du)$ is a mean-zero Gaussian random vector with covariance matrix equal to $h\int_0^1 K_h^2(u-v)\Sigma(u)\,du \to \nu_0\Sigma(v)$ as $h \to 0$. Hence, $h^{1/2}\int_0^1 K_h(u-v)\tilde{W}_I(du) \xrightarrow{\mathcal{D}} N(0, \nu_0\Sigma(v))$ as $h \to 0$, $nh \to \infty$. By the Slutsky theorem, $n^{-1/2}h^{1/2}U(v, \beta)$ converges weakly to $N(0, \nu_0\Sigma(v))$ as $nh^2 \to \infty$ and $nh^5 \to 0$.

Note that $U(\hat{\beta}(v)) - U(\beta(v)) = l_\beta''(v, \beta^*(v))(\hat{\beta}(v) - \beta(v))$, where $\beta^*(v)$ is on the line segment between $\hat{\beta}(v)$ and $\beta(v)$. By Condition A and the uniform consistency of $\hat{\beta}(v)$ on $v \in [a, b] \subset (0, 1)$, we have $n^{-1}l_\beta''(v, \beta^*(v)) = -\Sigma(v) + o_p(1)$, uniformly in $v \in [a, b]$ for $0 < \delta < 1/2$. Hence,

$$n^{1/2}h^{1/2}(\hat{\beta}(v) - \beta(v)) = -(l_\beta''(v, \beta^*(v))/n)^{-1}n^{-1/2}h^{1/2}U(\beta(v))$$

(17)
$$= (\Sigma(v))^{-1}n^{-1/2}h^{1/2}U(\beta(v)) + o_p(1),$$

uniformly in $v \in [a, b]$. It follows that $(nh)^{1/2}(\hat{\beta}(v) - \beta(v)) \xrightarrow{\mathcal{D}} N(0, \nu_0\Sigma(v)^{-1})$ as $nh^2 \to \infty$ and $nh^5 \to 0$. $\square$

PROOF OF THEOREM 4. From (16) and the first line of (17), we have, for $v \in [a, b]$,

$$\int_a^v n^{1/2}(\hat{\beta}(u) - \beta(u))\,du = -\int_a^v (\Sigma(u))^{-1}\int_0^1 K_h(x - u)\tilde{W}_I(dx)\,du + o_p(1).$$

Exchanging the order of integration and by the compact support of the kernel function $K(\cdot)$ on $[-1, 1]$, we have

$$\int_a^v n^{1/2}(\hat{\beta}(u) - \beta(u))\,du$$

$$= -\int_0^1 \left[\int_a^v (\Sigma(u))^{-1} K_h(x-u)\,du\right]\tilde{W}_I(dx) + o_p(1)$$

(18)
$$= -\int_{a+h}^{v-h} \left[\int_a^v (\Sigma(u))^{-1} K_h(x-u)\,du\right]\tilde{W}_I(dx)$$
$$- \int_{a-h}^{a+h} \left[\int_a^v (\Sigma(u))^{-1} K_h(x-u)\,du\right]\tilde{W}_I(dx)$$
$$- \int_{v-h}^{v+h} \left[\int_a^v (\Sigma(u))^{-1} K_h(x-u)\,du\right]\tilde{W}_I(dx) + o_p(1).$$

By Theorem 1, the process $\tilde{W}_I(x)$ converges weakly to a mean-zero Gaussian process with continuous paths. Under the assumption (A.4), $\int_a^v (\Sigma(u))^{-1} K_h(x-u)\,du$ has bounded variation and converges uniformly to $\Sigma(x)^{-1}$ for $x \in (a+h, v-h)$. By Lemma 2 of Gilbert, McKeague and Sun [5], the first term in (18) is $-\int_a^v (\Sigma(x))^{-1}\tilde{W}_I(dx) + o_p(1)$. Similar arguments lead to the second and the third terms in (18) to be $o_p(1)$. Hence

$$\int_a^v n^{1/2}(\hat{\beta}(u) - \beta(u))\,du = -\int_a^v (\Sigma(x))^{-1}\tilde{W}_I(dx) + o_p(1)$$
$$= -\tilde{W}_{\Sigma^{-1}}(v) + o_p(1),$$

which converges weakly to a $p$-dimensional mean-zero Gaussian martingale, $W_{\Sigma(v)^{-1}}(v)$, with continuous paths. The covariance matrix of $W_{\Sigma(v)^{-1}}(v)$ equals to $\mathrm{Cov}(W_{\Sigma^{-1}}(v)) = \int_a^v \Sigma(u)^{-1}\Sigma(u)\Sigma(u)^{-1}\,du = \int_a^v \Sigma(u)^{-1}\,du.$  $\square$

**Acknowledgments.** The authors gratefully acknowledge David Jobes and VaxGen Inc. for providing the HIV sequence data. The authors also thank the Associate Editor and two referees for their valuable comments.

Y. SUN  
DEPARTMENT OF MATHEMATICS AND STATISTICS  
UNIVERSITY OF NORTH CAROLINA AT CHARLOTTE  
9201 UNIVERSITY CITY BOULEVARD  
CHARLOTTE, NORTH CAROLINA 28223  
USA  
E-MAIL: yasun@uncc.edu  

P. B. GILBERT  
DEPARTMENT OF BIOSTATISTICS  
UNIVERSITY OF WASHINGTON  
AND  
FRED HUTCHINSON CANCER RESEARCH CENTER  
SEATTLE, WASHINGTON 98109  
USA  
E-MAIL: pgilbert@scharp.org

I. W. McKeague  
Department of Biostatistics  
Columbia University  
Mailman School of Public Health  
722 West 168th Street, 6th Floor  
New York, New York 10032  
USA  
E-mail: im2131@columbia.edu